\documentclass[11pt]{article}
\usepackage[english]{babel} 
\usepackage{a4wide} 
\usepackage{amssymb}
\usepackage{amsmath} 
\usepackage{enumitem} 
\usepackage{stmaryrd} 
\usepackage{bussproofs} 
\usepackage{amsthm} 
\usepackage[nottoc]{tocbibind} 
\usepackage{imakeidx} 
\usepackage{hyperref} 
\usepackage[capitalize]{cleveref} 
\crefformat{equation}{(#2#1#3)} 
\numberwithin{equation}{section}

\newtheorem{thm}{Theorem}[subsection]
\newtheorem{lem}[thm]{Lemma}
\newtheorem{cor}[thm]{Corollary}
\newtheorem{prop}[thm]{Proposition}

\theoremstyle{definition}
\newtheorem{defi}[thm]{Definition}

\newtheorem{rem}[thm]{Remark}
\newtheorem{conv}[thm]{Convention}

\newenvironment{prf}{\begingroup\setlength{\parskip}{0pt}\begin{proof}}{\end{proof}\endgroup}

\newcommand{\ha}{\mathsf{HA}}
\newcommand{\hahat}{{\widehat{\ha}}}
\newcommand{\ea}{{\mathsf{i}\s\textup{I}\Sigma_1}}
\newcommand{\s}{\hspace{1pt}}

\newcommand{\pa}{\mathsf{PA}}
\newcommand{\ax}{\mathsf{Ax}}
\newcommand{\bew}{\mathsf{Bew}}
\newcommand{\pf}{\mathsf{Prf}}
\newcommand{\iGLC}{\mathsf{iGLC}}
\newcommand{\ibew}{\mathsf{IBew}}
\renewcommand{\lg}{\ulcorner}
\newcommand{\rg}{\urcorner}
\newcommand{\nnil}{\mathsf{NNIL}}
\newcommand{\tnnil}{\mathsf{TNNIL}}
\newcommand{\sha}{\mathsf{sHA}}
\renewcommand{\r}{\sqsubset}
\renewcommand{\int}{\preceq}
\newcommand{\req}{\sqsubseteq}
\renewcommand{\th}{\mathsf{Th}}
\newcommand{\lval}{\llbracket}
\newcommand{\rval}{\rrbracket} 
\newcommand{\ruitje}{\hfill $\lozenge$} 
\newcommand{\wruitje}{\tag*{$\lozenge$}}

\setlist[itemize,enumerate]{itemsep=0pt,topsep=-6pt}

\begin{document}

\begin{center}
{\Large  Provability Logic and the Completeness Principle}
\end{center}

Albert Visser, Department of Philosophy and Religious Studies, Utrecht University\\
Jetze Zoethout, Department of Mathematics, Utrecht University

Version of \today

\begin{center}
\textbf{Abstract}
\end{center}\vspace{-1.5\baselineskip}
{\footnotesize
The logic $\iGLC$ is the intuitionistic version of L\"ob's Logic plus the completeness principle $A \to \Box A$.
In this paper, we prove an arithmetical completeness theorems for $\iGLC$ for theories equipped with two provability predicates $\Box$ and $\triangle$ that 
prove the schemes $A\to\triangle A$ and $\Box\triangle S\to\Box S$ for $S\in\Sigma_1$. We provide two salient instances of the theorem. In the first, $\Box$ is
fast provability and $\triangle$ is ordinary provability and, in the second, $\Box$ is ordinary provability and $\triangle$ is slow provability.

Using the second instance, we reprove a theorem previously obtained by 
M.\@ Ardeshir and S.\@ Mojtaba Mojtahedi \cite{ArdeshirMojtahedi} determining the $\Sigma_1$-provability logic of Heyting Arithmetic.

\textit{Keywords {\&} Phrases}. Provability Logic, Constructivism

\renewcommand{\baselinestretch}{0}
\tableofcontents
\renewcommand{\baselinestretch}{1}}

\newpage
\section{Introduction}

Around 1930, Kurt G\"odel proved his celebrated incompleteness theorems. While these results can be seen as the 
culmination of one era of logical research, they also cleared the way for several new fields within mathematical logic. 
An example of such a field is \emph{provability logic}, a topic that still occupies logicians today. Provability logic takes 
one of the main ingredients of G\"odel's theorems as its starting point. This ingredient is the formalization of the notion `formally provable in a certain 
arithmetical theory $T$' inside the language of arithmetic itself. Once this step has been taken, one may wonder what a theory $T$ is 
able to prove about its own notion of provability. This object, i.e.\@ what a theory $T$ can prove about its own notion of provability, is 
called the provability logic of $T$. Let us write, as we will below, `$\vdash_T A$' for `$A$ is formally provable in $T$', and `$\Box_TA$' for 
the arithmetical formula expressing that $A$ is formally provable in $T$. Then under some reasonable assumptions, the following turn out to hold:
\begin{itemize}
\item[(i)]	if $\vdash_T A$, then $\vdash_T\Box_TA$;
\item[(ii)]	$\vdash_T \Box_T(A\to B)\to (\Box_TA\to\Box_TB)$;
\item[(iii)]	$\vdash_T \Box_TA\to \Box_T\Box_TA$.
\end{itemize}
These are known nowadays as the \emph{Hilbert-Bernays-L\"ob derivability conditions}. Using another key idea from G\"odel's theorems, the 
Diagonalization Lemma, one can derive from these that $\vdash_T \Box_T(\Box_TA\to A)\to\Box_TA$, a result known as L\"ob's Theorem. In 1976, Robert Solovay proved that for the theory \emph{Peano Arithmetic}, the schemes (i)-(iii) and L\"ob's Theorem completely describe its provability logic \cite{Solovay}.

Provability logics are not monotone in their corresponding theories. That is, if $T$ is a theory extending another theory $U$, then it is not in general true that the provability logic of $T$ extends the provability logic of $U$. In light of this, it is all the more surprising that, in the classical case, provability logics are immensely stable. Solovay's proof can be modified to show that any $\Sigma_1$-sound theory interpreting Elementary Arithmetic has the same provability logic as Elementary Arithmetic. These include theories as strong as Zermelo-Fraenkel Set Theory (with or without the Axiom of Choice).

Elementary Arithmetic is a classical theory, which is why we made the caveat `in the classical case' above. In the intuitionistic case, the situation is completely different. Solovay's proof simply does not work for intuitionistic theories. This shows itself in the fact that the provability logic of \emph{Heyting Arithmetic}, the intuitionistic counterpart of Peano Arithmetic, contains principles that the provability logic of Peano Arithmetic does not share. These principles are somewhat exotic, and it is unknown what the provability logic of Heyting Arithmetic exactly is. In fact, as far as we are aware, there is presently only one intuitionistic theory for which a nontrivial provability logic is known, a result due to the first author (see \cite{completeness} and \cref{provlog_pa^ast}).

In Solovay's proof, the semantics of (classical) modal logic plays a major role. The larger part of the proof consists of embedding models for modal logic in a certain way into the theory $T$. These models are equipped with an accessibility relation. Solovay uses the predicate $\Box_T$ to represent this relation inside the theory $T$. One may try to give a Solovay-style proof by replacing the models for classical modal logic by models for intuitionistic modal logic. The difficulty about these models, however, is that they also possess an intuitionistic relation, in addition to the accessibility relation. The main question then becomes how we can deal with these two relations.

The main goal of this paper is to find interesting situations where we can give a Solovay-style embedding of a model for intuitionistic modal logic. 
We zoom in on models of the provability logic $\iGLC$. This logic is {\sf iGL}, the intuitionistic version of L\"ob's Logic, plus the Completeness Principle
$A\to \Box A$. This logic is, in a sense, the simplest interesting provability logic. Its models are pleasantly simple and hence good candidates for embeddability.

The Kripke models for $\iGLC$ have two accessibility relations.
To make the embedding work we use two notions of provability each one corresponding to one of the
accessibility relations. As a result we obtain arithmetical completeness results for $\iGLC$ in various theories and
for various interpretations of $\Box$.

Our Solovay-style embedding is presented in detail below. The embedding can be used to obtain a variety of results in provability logic. 
Among these is the determination of the \emph{$\Sigma_1$-provability logic of Heyting Arithmetic}, an object related to the ordinary 
provability logic of Heyting Arithmetic. This is not a new result. It was already obtained in 2014 by 
M.\@ Ardeshir and S.\@ Mojtaba Mojtahedi \cite{ArdeshirMojtahedi}, but the present work arrives at it in a different way. 
We stress, however, that our proof could not have been devised without the work from the paper \cite{ArdeshirMojtahedi}. 
First of all, it is of course easier to determine a provability logic if one already knows what it should be. Moreover, even though our proof is different, 
we do use some key ingredients from the paper \cite{ArdeshirMojtahedi}, most notably the $\tnnil$-algorithm.

Let us briefly outline the structure of the paper. First of all, in \cref{chap_prerequisites}, we discuss all the necessary prerequisite knowledge, and fix our notation. This section contains no essentially new results, but we do prove some results from the paper \cite{completeness} under weaker assumptions. For reasons of space, we will not spell out any specific G\"odel numberings or give an explicit definition of the predicate $\Box_T$. Therefore, it will be useful to have some prior experience with G\"odel's incompleteness theorems and with provability logic (in the classical case) when reading this paper. A reader that is already familiar with (some of) the concepts discussed in \cref{chap_prerequisites} may want to read (a portion of) this chapter only superficially, and refer back to it if necessary. In \cref{chap_CT}, we present our Solovay-style embedding, and formulate our central
 completeness theorem. This theorem will be stated in an abstract way that does not yet mention any specific theories or provability predicates. In \cref{chap_applications}, we will present several applications of our completeness theorem, among which the determination of the $\Sigma_1$-provability logic of Heyting Arithmetic.

An earlier version of the current paper was submitted by the second author as a thesis for the MA degree in philosophy at Utrecht University. This thesis was supervised by the first author. We would like to thank Rosalie Iemhoff, Lev Beklemishev and Sven Bosman for commenting on drafts of this thesis and the current paper. We are grateful to Mojtaba Mojtahedi for his comments on the preprint version.

\newpage\section{Prerequisites}\label{chap_prerequisites}

In this section, we develop some notation and theory that will be used in the later parts of the paper. First, in \cref{sec_arithmetic_provability}, we fix some basic notions about arithmetical theories and provability predicates. Then, in \cref{sec_CP}, we discuss the $T$-translation, which will lead to theories that prove their own completeness. In \cref{sec_fast_slow}, we turn our attention to two nonstandard notions of provability, called fast and slow provability. Finally, in \cref{sec_IML}, we develop some intuitionistic (propositional) modal logic.

\subsection{Arithmetic and Provability}\label{sec_arithmetic_provability}

All the theories we shall consider will be theories for intuitionistic predicate logic with equality. As our proof system, we pick natural deduction with equality.\index[term]{Natural deduction with equality} An \emph{axiom}\index[term]{Axiom} will be viewed as a special case of an inference rule, namely as an inference rule whose premiss set is empty. For equality, we have the axiom $x=x$, and an inference rule involving substitution. The language in which our theories will be fomulated will be the \emph{language of arithmetic} $\mathcal{L}=\{0,\mathsf{S},+,\times\}$.\index[term]{Language!of arithmetic}\index[symb]{al@$\mathcal{L}$}\index[symb]{as@$\mathsf{S}$} Here 0 is a constant symbol, $\mathsf{S}$ is a unary function symbol and $+$ and $\times$ are binary function symbols. For each $n\in\mathbb{N}$, we can define the $\mathcal{L}$-term $\mathsf{S}\ldots \mathsf{S}0$, where the $\mathsf{S}$ occurs exactly $n$ times. This term is called the \emph{numeral}\index[term]{Numeral} of $n$, and we denote it just by $\overline{n}$. For terms $s$ and $t$, we define $s\leq t$ as $\exists x\s (s+x=t)$ and $s<t$ as $\exists x\s (s+\mathsf{S}x=t)$.\index[symb]{asleqt@$s\leq t$}\index[symb]{asleqt1@$s<t$} Here $x$ should not occur in $s$ or $t$, of course. We notice that the language $\mathcal{L}$ has a straightforward interpretation in the natural numbers, yielding the \emph{standard model}\index[term]{Standard model} $\mathbb{N}$. We introduce two special classes of formulae.
\begin{defi}\index[symb]{adelta0@$\Delta_0$}\index[symb]{asigma1@$\Sigma_1$}
\begin{itemize}
\item[(i)]	The set of \emph{$\Delta_0$-formulae} is defined by recursion, as follows: \vspace{.25\baselineskip}
\begin{itemize}
\item[(a)]	all atomic $\mathcal{L}$-formulae are $\Delta_0$-formulae;
\item[(b)]	the set of $\Delta_0$-formulae is closed under conjunction, disjunction and implication;
\item[(c)]	if $A$ is a $\Delta_0$-formula, and $t$ is an $\mathcal{L}$-term not containing the variable $x$, then the formulae $\exists x\s(x<t\wedge A)$ and $\forall x\s(x<t\to A)$ are also $\Delta_0$-formulae. \vspace{.25\baselineskip}
\end{itemize}
We write $A\in\Delta_0$ if $A$ is a $\Delta_0$-formula.
\item[(ii)]	The set of \emph{$\Sigma_1$-formulae} consists of all $\mathcal{L}$-formulae of the form $\exists x\s A$, where $A\in\Delta_0$. We write $S\in\Sigma_1$ if $S$ is a $\Sigma_1$-formula. \ruitje
\end{itemize}
\end{defi}
To each $\mathcal{L}$-expression $\alpha$ (which can be a term, a formula or a sequence of formulae), we assign a \emph{G\"odel number}\index[term]{Godel@G\"odel number} $\lg\alpha\rg$ in some reasonable way. More precisely, we require that elementary syntactic operations concerning $\mathcal{L}$ are elementary functions (to be defined below) in terms of their G\"odel numbers.

\begin{defi}\index[term]{Theory}\index[symb]{aaxt@$\ax_T$}\index[symb]{atht@$\th(T)$}
A \emph{theory} $T$ will be a pair $(\th(T),\ax_T)$, where $\ax_T$ is a $\Sigma_1$-formula in one free variable, and $\th(T)$ is precisely the set of $\mathcal{L}$-formulae derivable from the axiom set 
\[
\{ A\mid A\mbox{ an }\mathcal{L}\mbox{-formula, }\mathbb{N}\models \ax_T(\lg A\rg)\}. \wruitje
\]
\end{defi}

In other words, a theory is a set of $\mathcal{L}$-formulae that is closed under derivability in intuitionistic predicate logic with equality, together with a $\Sigma_1$-formula that defines an axiom set for the theory in the standard model. Usually, we will define a theory by giving its axioms, understanding that their is some natural $\Sigma_1$-formulation in $\mathcal{L}$ for axiomhood. For a set of $\mathcal{L}$-formulae $\Gamma$ and an $\mathcal{L}$-formula $A$, we write $\Gamma\vdash_T A$\index[symb]{avdasht@$\vdash_T$!for a theory} to indicate that $A$ is provable using open assumptions from $\Gamma$ and the axioms of $T$. Notice that $\vdash_T A$ just means $A\in \th(T)$. Now we define three theories that will be of great interest to us. 
\begin{defi}\label{def_ha}
\begin{itemize}\index[symb]{aha@$\ha$}\index[symb]{aiisigma1@$\ea$}\index[symb]{apa@$\pa$}\index[term]{Induction axiom}\index[term]{Heyting Arithmetic}\index[term]{Peano Arithmetic}\index[term]{Law of the Excluded Middle}
\item[(i)]	The theory $\ea$, called \emph{intuitionistic $\Sigma_1$-induction}, has the axioms

\begin{tabular}{ll}
$\neg(\mathsf{S}x=0)$ \hspace{3.5cm} & $\mathsf{S}x=\mathsf{S}y\to x=y$\\
$x+0=0$ & $x+\mathsf{S}y = \mathsf{S}(x+y)$\\
$x\times 0=0$ & $x\times \mathsf{S}y = x\times y + x$
\end{tabular}

and, for each $\mathcal{L}$-formula $S\in\Sigma_1$, the \emph{induction axiom} 
\[S[0/x] \wedge \forall x\s (S\to S[\mathsf{S}x/x]) \to \forall x\s S.\]
\item[(ii)]	The theory $\ha$, called \emph{Heyting arithmetic}, has all the axioms of $\ea$, together with \emph{full induction}: for each $\mathcal{L}$-formula $A$, we have the axiom
\[A[0/x] \wedge \forall x\s (A\to A[\mathsf{S}x/x]) \to \forall x\s A.\]
\item[(iii)]	The theory $\pa$, called \emph{Peano arithmetic}, has all the axioms of $\ha$, together with the \emph{Law of the Excluded Middle}: $A\vee\neg A$, where $A$ is an $\mathcal{L}$-formula. \ruitje
\end{itemize}
\end{defi}
Even though the axiom set we presented for $\ea$ is infinite, the theory $\ea$ is actually finitely axiomatizable. This is because the induction scheme for $\Sigma_1$-formulae follows the single induction axiom for the $\Sigma_1$-truth predicate, which is itself a $\Sigma_1$-formula. The finite axiomatizability of $\ea$ is verifiable in $\ea$ itself. It is also well-known that $\ea$, and hence any theory extending it, is \emph{$\Sigma_1$-complete}.\index[term]{sigma1completeness@$\Sigma_1$-completeness} That is, every $\Sigma_1$-sentence true in the standard model can be proven inside $\ea$. 

Even though $\ea$ is an intuitionistic theory, we do have the following result, which can be shown using induction on complexity.
\begin{prop}\label{Delta0_stuff}
If $A\in\Delta_0$ is a formula, then $\vdash_\ea A\vee\neg A$.
\end{prop}
We also have the following result, that we shall not prove.
\begin{prop}\label{representing_el}
Let $F\colon \mathbb{N}^k\to\mathbb{N}$ be a primitive recursive function. Then there exists a $\Sigma_1$-formula $A_F(\vec{x},y)$ satisfying:
\begin{itemize}
\item[\textup{(}i\textup{)}]	$\vdash_\ea A_F(\vec{n},F(\vec{n}))$ for all $\vec{n}\in\mathbb{N}^k$;
\item[\textup{(}ii\textup{)}]	$\vdash_\ea \exists y\forall z\s (A_F(\vec{x},z)\leftrightarrow y=z)$.
\end{itemize}
Moreover, this formula can be chosen in such a way that the definition of $F$ as a primitive recursive function is verifiable in $\ea$.
\end{prop}
We have a primitive recursive function $\text{Subst}\colon\mathbb{N}^2\to\mathbb{N}$ that is defined as follows.\index[symb]{asubst@Subst} If $a$ is the G\"odel number of some formula $A(v)$ in one free variable $v$, then $\text{Subst}(a,b)=\lg A(b)\rg$; otherwise, $\text{Subst}(a,b)=0$. We can represent this function in $\ea$ using \cref{representing_el}. If $A(v)$ is a formula with one free variable, we will write $\lg A(\tilde{x})\rg$ for $\text{Subst}(\lg A(v)\rg,x)$, which makes sense when working in a theory extending $\ea$. We apply similar conventions for multiple free variables.\index[symb]{axtilde@$\tilde{x}$} We will need the following famous result, that we will not prove.
\begin{thm}[Diagonalization Lemma]\index[term]{Diagonalization Lemma}
Suppose $A(\vec{x},y)$ is an $\mathcal{L}$-formula. Then there exists an $\mathcal{L}$-formula $B(\vec{x})$ such that $\vdash_\ea B(\vec{x}) \leftrightarrow A(\vec{x},\lg B(\vec{x})\rg)$.
\end{thm}
Now suppose we have a theory $T$. Using the $\Sigma_1$-formula $\ax_T$, we can construct a $\Sigma_1$-formula $\bew_T(x)$\index[symb]{abewt@$\bew_T$} that expresses `$x$ is the G\"odel number of some formula $A$ such that $\vdash_T A$' in a natural way. We can write $\bew(x)$ as $\exists y\s \pf_T(y,x)$\index[symb]{aprft@$\pf_T$} for some $\Delta_0$-formula $\pf_T$. We think of $\pf(y,x)$ as expressing the fact that $y$ codes a $T$-proof of the formula that has $x$ as its G\"odel number. For a formula $A=A(x_1, \ldots, x_n)$, we write $\Box_T A$\index[symb]{aboxt@$\Box_T$} for $\bew_T(\lg A(\tilde{x}_1,\ldots,\tilde{x}_n)\rg)$. In particular, $\Box_T A$ has the same free variables as $A$. Now we can define certain relations between theories.
\begin{defi}\label{comparing_theories}
Let $U$ and $T$ be theories. We write:
\begin{itemize}\index[symb]{aut1@$U\subseteq T$} \index[symb]{aut2@$U=T$} \index[symb]{aut3@$U\leq T$} \index[symb]{aut4@$U\equiv T$}
\item[(i)]	$U\subseteq T$ if $\th(U)\subseteq \th(T)$;
\item[(ii)]	$U=T$ if $\th(U)=\th(T)$;
\item[(iii)]	$U\leq T$ if $\vdash_\ea \bew_U(x)\to\bew_T(x)$;
\item[(iv)]	$U\equiv T$ if $\vdash_\ea \bew_U(x)\leftrightarrow\bew_T(x)$. \ruitje
\end{itemize}
\end{defi}
We emphasize that, then we write $U=T$, we do not mean an equality of the pairs $(\th(U),\ax_U)$ and $(\th(T),\ax_T)$, but only an equality of the first coordinate. Since $\ea$ is sound, we see that $U\leq T$ implies that $U\subseteq T$. We also notice that, if $U$ and $T$ are theories such that $\vdash_\ea\ax_U(x)\to \ax_T(x)$, then $U\leq T$ clearly holds. However, this requirement is not necessary: it can also be the case that every $U$-proof can (verifiably in $\ea$) be transformed into a $T$-proof without the one axiom set being contained in the other. Before we can develop more theory, we need to restrict our investigation to theories that, verifiably in $\ea$, can perform a minimal amount of arithmetic.
\begin{conv}\label{convention_on_T}\index[term]{Convention on theories}
All the theories $T$ we shall consider, will satisfy $\ea\leq T$. \ruitje
\end{conv}
\begin{rem}
As \cref{comparing_theories} and \cref{convention_on_T} make clear, $\ea$ will serve as our base theory. Most, but not all, of the following goes through for slightly weaker base theories, such as (intuitionistic) Elementary Arithmetic extended with $\Sigma_1$-collection. We have chosen $\ea$ because this yields the most simple and uniform exposition of the material below. \ruitje
\end{rem}
Notice that this clearly holds for the three theories from \cref{def_ha}. With this requirement in place, we can state some basic properties of $\Box_T$, that we will not prove.
\begin{prop}\label{basic_facts_BoxT}
Let $T$ be a theory and let $A$, $B$ and $S$ be $\mathcal{L}$-formulae. Then we have:
\begin{itemize}\index[term]{sigma1completeness@$\Sigma_1$-completeness!Formalized $\Sigma_1$-completeness} \index[term]{Lobt@L\"ob's Theorem}\index[term]{Lobp@L\"ob's Principle}
\item[\textup{(}i\textup{)}]	$\vdash_T A$ if and only if $\mathbb{N}\models \Box_TA$, if and only if $\vdash_\ea \Box_TA$;
\item[\textup{(}ii\textup{)}]	$\vdash_\ea \Box_T(A\to B)\to (\Box_TA\to\Box_TB)$;
\item[\textup{(}iii\textup{)}]	$\vdash_\ea \Box_TA\to\Box_T\Box_TA$;
\item[\textup{(}iv\textup{)}]	\textup{(Formalized $\Sigma_1$-completeness)} if $S\in\Sigma_1$, then $\vdash_\ea S\to \Box_T S$;
\item[\textup{(}v\textup{)}]	\textup{(L\"ob's Principle)} if $U$ is a theory such that $U\subseteq T$ and $\vdash_U \Box_T A\to A$, then $\vdash_U A$;
\item[\textup{(}vi\textup{)}]	\textup{(L\"ob's Theorem)} $\vdash_\ea \Box_T(\Box_T A\to A)\to\Box_TA$.
\end{itemize}
Moreover, \textup{(}ii\textup{)}, \textup{(}iii\textup{)}, \textup{(}iv\textup{)} and \textup{(}vi\textup{)} are verifiable in $\ea$.
\end{prop}
We remark that for (iii)-(vi), we need \cref{convention_on_T}. In the next section, we will need the following facts.
\begin{prop}\label{inclusion_and_Boxrule}
Let $U$ and $T$ be theories.
\begin{itemize}
\item[\textup{(}i\textup{)}]	If $U\subseteq T$, then $\vdash_U A$ implies $\vdash_U \Box_TA$ for all $\mathcal{L}$-formulae $A$.
\item[\textup{(}ii\textup{)}]	If $U\leq T$, then $\vdash_\ea \bew_U(x)\to\Box_U\bew_T(x)$. In particular, $\vdash_\ea \Box_UA\to\Box_U\Box_TA$ for all $\mathcal{L}$-formulae $A$.
\end{itemize}
\end{prop}
\begin{prf}
(i) If $\vdash_UA$, then also $\vdash_TA$, so $\vdash_{\ea} \Box_TA$. Since $\ea\subseteq U$, we also get $\vdash_U\Box_TA$.

(ii) Since $U\leq T$, we have $\vdash_\ea \bew_U(x)\to \bew_T(x)$. Since $\ea\subseteq U$ (by \cref{convention_on_T}), it follows from (i) that $\vdash_\ea \Box_U\bew_U(x)\to\Box_U\bew_T(x)$. We also have $\vdash_\ea\bew_U(x)\to\Box_U\bew_U(x)$ by formalized $\Sigma_1$-completeness, and now the result follows.
\end{prf}
For future use, we state the following definition.
\begin{defi}\index[term]{Provability predicate}
Let $T$ be a theory and let $P(x)$ be a $\Sigma_1$-formula in one free variable. For an $\mathcal{L}$-sentence $A$, we write $\Box A$\index[symb]{abox@$\Box$!for a provability predicate} for $P(\lg A\rg)$. We say that $P$ is a \emph{provability predicate} for $T$ if the following hold for all $\mathcal{L}$-sentences $A$, $B$ and $S$:
\begin{itemize}
\item[\textup{(}i\textup{)}]	if $\vdash_TA$, then $\mathbb{N}\models \Box A$;
\item[\textup{(}ii\textup{)}]	$\vdash_\ea \Box(A\to B)\to (\Box A\to\Box B)$;
\item[\textup{(}iii\textup{)}]	if $S\in\Sigma_1$, then $\vdash_\ea S\to\Box S$. \ruitje
\end{itemize}
\end{defi}
Observe that the above definition has the following monotonicity property: if $T'\subseteq T$ are theories and $P$ is a provability predicate for $T$, then $P$ is also a provability predicate for $T'$. Using \cref{basic_facts_BoxT}, we see that $\bew_T$ is always a provability predicate for $T$. In fact, any provability predicate for $T$ has properties similar to those of $\bew_T$, as the following result, whose proof is standard, shows.
\begin{prop}\label{properties_pp}
Let $P$ be provability predicate for a theory $T$. For $\mathcal{L}$-sentences $A$, write $\Box A$ for $P(\lg A\rg)$. Then for all $\mathcal{L}$-sentences $A$, we have:
\begin{itemize}
\item[\textup{(}i\textup{)}]	$\vdash_\ea \Box A\to\Box\Box A$;
\item[\textup{(}ii\textup{)}]	if $U$ is a theory such that $U\subseteq T$ and $\vdash_U \Box A\to A$, then $\vdash_U A$;
\item[\textup{(}iii\textup{)}]	$\vdash_\ea \Box(\Box A\to A)\to\Box A$.
\end{itemize}
\end{prop}

\subsection{The Completeness Principle}\label{sec_CP}

In this section, we introduce the $T$-translation, that will allow us to define theories that prove their own completeness. All results in this section are from the paper \textit{On the Completeness Principle} \cite{completeness}, but we have formulated some of them under weaker conditions.
\begin{defi}\index[term]{T-translation@$T$-translation}\index[symb]{aat5@$A^T$}
Let $T$ be a theory. We define the \emph{$T$-translation} $(\cdot)^T$ from the set of $\mathcal{L}$-formulae to itself by recursion. For all $\mathcal{L}$-terms $s$ and $t$ and $\mathcal{L}$-formulae $A$ and $B$, we set:
\begin{itemize}
\item[(i)]	$(s=t)^T$ is $s=t$ and $\perp^T$ is $\perp$;
\item[(ii)]	$(A\circ B)^T$ is $A^T\circ B^T$ for $\circ\in\{\wedge,\vee\}$;
\item[(iii)]	$(A\to B)^T$ is $(A^T\to B^T)\wedge\Box_T(A^T\to B^T)$;
\item[(iv)]	$(\exists x\s A)^T$ is $\exists x\s A^T$;
\item[(v)]	$(\forall x\s A)^T$ is $\forall x\s A^T\wedge\Box_T(\forall x\s A^T)$. \ruitje
\end{itemize}
\end{defi}
Based on the $T$-translation, we can construct new theories out of existing ones.
\begin{defi}\index[symb]{aut@$U^T$}\index[symb]{austar@$U^\ast$}
Let $U$ and $T$ be theories. We define the theory $U^T$ as the theory whose axioms are those of $\ea$, together with the set $\{A\mid\ \vdash_U A^T\}$. For a theory $U$, we write $U^\ast$ for $U^U$. \ruitje
\end{defi}
We make some remarks on how $\ax_{U^T}$ can be defined. The function $(\cdot)^T\colon\mathbb{N}\to\mathbb{N}$ that satisfies $x^T = \lg A^T\rg$ if $x$ is the G\"odel number of an $\mathcal{L}$-formula $A$, and $x^T=0$ otherwise, is primitive recursive. This means that we can represent this function in $\ea$ using \cref{representing_el}. Now we define $\ax_{U^T}(x)$ as $\ax_\ea(x)\vee (\text{Form}(x)\wedge\bew_U(x^T))$, where $\text{Form}(x)\in\Sigma_1$ naturally expresses the fact that $x$ is the G\"odel number of an $\mathcal{L}$-formula.\index[symb]{aform@Form} We study the relation between provability in $U^T$ and provability in $U$ through the following lemmata. Our first lemma is the \emph{raison d'\^etre} for the $T$-translation.
\begin{lem}\label{completeness_for_translations}
For all $\mathcal{L}$-formulae $A$, we have $\vdash_\ea A^T\to\Box_T A^T$. Moreover, this is verifiable in $\ea$.
\end{lem}
\begin{prf}
We proceed by induction on the complexity of $A$.
\begin{itemize}
\item[At]	If $A$ is atomic, then $A^T=A$ and the claim follows from \cref{basic_facts_BoxT}(iv) since $A$ is a $\Sigma_1$-formula.
\item[$\wedge$]	Suppose $A=B\wedge C$ and the claim holds for $B$ and $C$. Then $A^T$ is $B^T\wedge C^T$, and we have
\[\vdash_\ea B^T\wedge C^T\to \Box_TB^T\wedge\Box_TC^T\to \Box_T(B^T\wedge C^T),\]
as desired.
\item[$\vee$]	Suppose $A$ is $B\vee C$ and the claim holds for $B$ and $C$. Then $A^T$ is $B^T\vee C^T$, and we have $\vdash_\ea B^T\to \Box_TB^T \to \Box_T(B^T\vee C^T)$ and $\vdash_\ea C^T\to \Box_TC^T \to \Box_T(B^T\vee C^T)$, which together yield $\vdash_\ea B^T\vee C^T\to \Box_T(B^T\vee C^T)$, as desired.
\item[$\to$]	Suppose $A$ is $B\to C$ and the claim holds for $B$ and $C$. Then the formula $A^T$ is equal to $(B^T\to C^T)\wedge\Box_T(B^T\to C^T)$, and we have
\begin{align*}
\vdash_\ea (B^T\to C^T)\wedge\Box_T(B^T\to C^T) &\to \Box_T(B^T\to C^T)\\
&\to \Box_T(B^T\to C^T)\wedge\Box_T\Box_T(B^T\to C^T)\\
&\to \Box_T((B^T\to C^T)\wedge\Box_T(B^T\to C^T)),
\end{align*}
as desired.
\item[$\exists$]	Suppose $A$ is $\exists x\s B$ and the claim holds for $B$. Then $A^T$ is $\exists x\s B^T$. It is provable in intuitionistic predicate logic that $B^T\to \exists x\s B^T$, so we also have $\vdash_\ea \Box_TB^T\to\Box_T(\exists x\s B^T)$. We get $\vdash_\ea \exists x\s B^T \to \exists x\s\Box_T B^T \to \Box_T(\exists x\s B^T)$, as desired.
\item[$\forall$]	Suppose $A$ is $\forall x\s B$ and the claim holds for $B$. Then $A^T$ is $\forall x\s B^T\wedge\Box_T(\forall x\s B^T)$, and we have
\begin{align*}
\vdash_\ea \forall x\s B^T\wedge\Box_T(\forall x\s B^T) &\to \Box_T(\forall x\s B^T)\\
&\to \Box_T(\forall x\s B^T)\wedge \Box_T\Box_T(\forall x\s B^T)\\
&\to \Box_T(\forall x\s B^T\wedge\Box_T(\forall x\s B^T)),
\end{align*}
as desired.\vspace{.5\baselineskip}
\end{itemize}
For the second statement, we should carry out this induction inside $\ea$. One should notice that now we need that clauses (ii)-(iv) from \cref{basic_facts_BoxT} are verifiable in $\ea$.
\end{prf}
Next we show that, up to equivalence, $\Sigma_1$-formulae are invariant under the $T$-translation.
\begin{lem}\label{Sigma_and_T}
If $T$ is a theory and $S$ is a $\Sigma_1$-formula, then $\vdash_\ea S\leftrightarrow S^T$. Moreover, this is verifiable in $\ea$.
\end{lem}
\begin{prf}
If $A$ is a $\Delta_0$-formula, then $\vdash_\ea A\to\Box_T A$, which means that $\vdash_\ea A\wedge\Box_T A\leftrightarrow A$. Using this observation, we can show that $\vdash_\ea A\leftrightarrow A^T$ for all $A\in\Delta_0$ using a straightforward induction on the complexity of $A$. Finally, if $S\in\Sigma_1$, then write $S$ as $\exists x\s A$ with $A\in\Delta_0$. Then $S^T$ is $\exists x\s A^T$, and we see that $\vdash_\ea \exists x\s A\leftrightarrow\exists x\s A^T$, as desired.

For the second statement, we formalize the above in $\ea$.
\end{prf}
Finally, we collect some technical facts about the interaction between the $T$-translation and substitution.
\begin{lem} \label{technical_stuff}
Let $A$ be a formula, let $x$ be a variable, and let $s$ be a term. Then:
\begin{itemize}
\item[\textup{(}i\textup{)}]	$A$ and $A^T$ have the same free variables;
\item[\textup{(}ii\textup{)}]	$s$ is free for $x$ in $A$ if and only if $s$ is free for $x$ in $A^T$;
\item[\textup{(}iii\textup{)}]	if $s$ is free for $x$ in $A$, then $\vdash_\ea (A^T)[s/x] \leftrightarrow (A[s/x])^T$.
\end{itemize}
Moreover, these are all verifiable in $\ea$.
\end{lem}
\begin{prf} All three statements can be proven by an easy induction on the complexity of $A$. For the induction steps for implication and universal quantification in statement (iii), one should observe that, verifiably in $\ea$, we have $\vdash_\ea (\Box_TA)[s/x]\leftrightarrow \Box_T(A[s/x])$ for all $\mathcal{L}$-terms $s$ and $\mathcal{L}$-formulae $A$.
\end{prf}
Using these lemmata, we can prove the following crucial result.
\begin{thm}\label{prov_in_U_and_U^T}
Let $U$ and $T$ be theories such that $\vdash_U B$ implies $\vdash_U \Box_T B$ for all $\mathcal{L}$-formulae $B$. For a set of $\mathcal{L}$-formulae $\Gamma$, write $\Gamma^T=\{B^T\mid B\in\Gamma\}$. Then for all $\mathcal{L}$-formulae $A$, we have $\Gamma\vdash_{U^T} A$ if and only if $\Gamma^T\vdash_U A^T$.
\end{thm}
\begin{rem}
\begin{itemize}
\item[(i)]	By \cref{inclusion_and_Boxrule}(i), the conditions on $U$ and $T$ apply in particular when $U\subseteq T$. We formulate this theorem (and \cref{formalized_prov_in_U_and_U^T} below) in such a strong way in order to obtain \cref{prov_in_hahat} below.
\item[(ii)]	We warn the reader that, under these conditions on $U$ and $T$, we cannot necessarily verify the result `$\Gamma\vdash_{U^T} A$ if and only if $\Gamma^T\vdash_U A^T$' inside $\ea$; see \cref{formalized_prov_in_U_and_U^T} below. \ruitje
\end{itemize}
\end{rem}
\begin{prf}[Proof of \cref{prov_in_U_and_U^T}]
Suppose that $\Gamma^T\vdash_U A^T$. Then there exist $n\geq 0$ and $C_1, \cdots, C_n\in\Gamma$ such that $\vdash_U C^T_0\wedge\ldots\wedge C^T_n\to A^T$. Then we also have $\vdash_U(C_0\wedge\ldots\wedge C_n)^T\to A^T$, and by our assumption, we also get $\vdash_U\Box_T((C_0\wedge\ldots\wedge C_n)^T\to A^T)$. So $\vdash_U (C_0\wedge\ldots\wedge C_n\to A)^T$, and therefore we get $\vdash_{U^T} C_0\wedge\ldots\wedge C_n\to A$. Finally, this clearly yields that $\Gamma\vdash_{U^T} A$.\\

For the converse direction, we proceed by induction on the proof tree for $\Gamma\vdash_{U^T} A$. Before we start, we notice the following: if $\vdash_U B\to C$ for certain $\mathcal{L}$-formulae $B$ and $C$, then by our assumption, $\vdash_U \Box_T(B\to C)$. We also have $\vdash_U \Box_T(B\to C)\to(\Box_TB\to\Box_TC)$, so we get $\vdash_U\Box_TB\to\Box_TC$. We also note: if $\vdash_\ea B$, then $\vdash_U B$, whence $\vdash_U\Box_TB$.\\

First, suppose that $A$ is an axiom of $U^T$. That is, we suppose that $A$ is an axiom of $\ea$ or that $\vdash_U A^T$. In the latter case, we are done. So suppose that $A$ is an axiom of $\ea$. We need to show that $\vdash_U A^T$. If $A$ is the axiom $x=x$ or a basic axiom of $\ea$, then $A\in\Delta_0$, so by \cref{Sigma_and_T}, we have $\vdash_\ea A\leftrightarrow A^T$. Since $\vdash_\ea A$, we also get $\vdash_\ea A^T$, and in particular, $\vdash_U A^T$. It remains to prove the claim for the case where $A$ is an induction axiom, say $S[0/x]\wedge\forall x\s(S\to S[\mathsf{S}x/x])\to \forall x\s S$ with $S\in\Sigma_1$. First of all, we notice that
\begin{align}
\vdash_\ea (\forall x\s (S\to S[\mathsf{S}x/x]))^T &\leftrightarrow \forall x\s (S\to S[\mathsf{S}x/x])^T\wedge\Box_T(\forall x\s  (S\to S[\mathsf{S}x/x])^T)\nonumber\\
&\to  \forall x\s (S^T\to (S[\mathsf{S}x/x])^T)\wedge\Box_T(\forall x\s (S^T\to (S[\mathsf{S}x/x])^T))\nonumber\\
&\leftrightarrow  \forall x\s (S^T\to (S^T)[\mathsf{S}x/x])\wedge\Box_T(\forall x\s (S^T\to (S^T)[\mathsf{S}x/x])).\label{translation_induction_axiom}
\end{align}
Furthermore, we know that $\vdash_\ea (S[0/x])^T\leftrightarrow (S^T)[0/x]$ and that $(\forall x\s S)^T$ is the formula $\forall x\s S^T\wedge\Box_T(\forall x\s S^T)$. Define the formulae
\begin{align*}
C &:\leftrightarrow (S^T)[0/x]\wedge\forall x\s(S^T\to (S^T)[\mathsf{S}x/x])\wedge\Box_T(\forall x\s(S^T\to (S^T)[\mathsf{S}x/x]))\to \forall x \s S^T\wedge\Box_T(\forall x\s S^T),\\
D &:\leftrightarrow (S[0/x])^T\wedge (\forall x\s (S\to S[\mathsf{S}x/x]))^T \to (\forall x\s S)^T.
\end{align*}
Then it follows from \cref{translation_induction_axiom} that $\vdash_\ea C\to D$. Now we also get $\vdash_U C\to D$ and hence $\vdash_U\Box_TC\to \Box_TD$. Since $A^T$ is the formula $D\wedge\Box_TD$, we see that $\vdash_U C\wedge\Box_T C\to A^T$. So it suffices to show that $\vdash_UC$.

Since $S\in\Sigma_1$, we have $\vdash_\ea S\leftrightarrow S^T$. This means that the induction axiom for $S^T$ is provable in $\ea$, hence in $U$:
\begin{align}\label{induction_axiom}
\vdash_U (S^T)[0/x]\wedge\forall x\s (S^T\to (S^T)[\mathsf{S}x/x])\to \forall x\s S^T.
\end{align}
Now it follows that
\begin{align}\label{boxed_induction_axiom}
\vdash_U \Box_T((S^T)[0/x])\wedge\Box_T(\forall x\s (S^T\to (S^T)[\mathsf{S}x/x]))\to \Box_T(\forall x\s S^T).
\end{align}
Finally, since $\vdash_\ea (S^T)[0/x] \leftrightarrow (S[0/x])^T$, we can use \cref{completeness_for_translations} to see that
\begin{align}\label{induction_base}
\vdash_\ea (S^T)[0/x]\to (S[0/x])^T \to \Box_T((S[0/x])^T)\to \Box_T((S^T)[0/x]).
\end{align}
From \cref{induction_axiom}, \cref{boxed_induction_axiom} and \cref{induction_base}, we may deduce that $C$ is indeed provable in $U$, as desired.\\

Now we treat the rules of inference. Since the $T$-translation commutes with conjunction, disjunction and existential quantification, the induction steps for rules of inference for these operators are trivial. It remains to check the rules for implication and universal quantification, and the substitution rule.
\begin{itemize}
\item[$\to$E]	Suppose that $\Gamma^T\vdash_U (B\to C)^T$ and $\Gamma^T\vdash_U B^T$. We need to show that $\Gamma^T\vdash_U C^T$. But this is obvious since $\vdash_\ea (B\to C)^T\to(B^T\to C^T)$.
\item[$\to$I]	Suppose that $\Gamma^T, B^T\vdash_U C^T$. We need to show that $\Gamma^T\vdash_U (B\to C)^T$. We certainly have $\Gamma^T\vdash_U B^T\to C^T$. But then we also have $\Box_T\Gamma^T\vdash_U \Box_T(B^T\to C^T)$, where $\Box_T\Gamma^T = \{\Box_TD^T\mid D\in\Gamma\}$. Since $\vdash_\ea D^T\to\Box_T D^T$ for all $D\in\Gamma$, we get $\Gamma^T\vdash_U\Box_T(B^T\to C^T)$. Combining our results, we find 
\[\Gamma^T\vdash_U (B^T\to C^T)\wedge\Box_T(B^T\to C^T),\]
as desired.
\item[$\forall$E]	Suppose that $\Gamma^T\vdash_U (\forall x\s B)^T$. We need to show that $\Gamma^T\vdash_U (B[s/x])^T$. Since $\vdash_\ea (\forall x\s B)^T\to \forall x\s B^T$, we see that $\Gamma^T\vdash_U (B^T)[s/x]$. Since we also know that $\vdash_\ea (B^T)[s/x]\leftrightarrow (B[s/x])^T$, we get $\Gamma^T\vdash_U (B[s/x])^T$, as desired.
\item[$\forall$I]	Suppose that $\Gamma^T\vdash_U B^T$, where the variable $x$ does not occur anywhere in $\Gamma$. We need to show that $\Gamma^T\vdash_U (\forall x\s B)^T$. First of all, we certainly have $\Gamma^T\vdash_U\forall x\s B^T$, since $x$ does not occur free anywhere in $\Gamma^T$. By applying the same reasoning as in the $\to$I-case, we find $\Gamma^T\vdash_U\Box_T(\forall x\s B^T)$. We conclude that $\Gamma^T\vdash_U\forall x\s B^T\wedge\Box_T(\forall x\s B^T)$, as desired.
\item[Subst]	Suppose that $\Gamma^T\vdash_U (B[s/x])^T$ and $\Gamma^T\vdash_U (s=t)^T$. We need to show that $\Gamma^T\vdash_U (B[t/x])^T$. We have $\Gamma^T\vdash_U s=t$ and by \cref{technical_stuff}(iii), we get $\Gamma^T\vdash_U (B^T)[s/x]$. This yields $\Gamma^T\vdash_U (B^T)[t/x]$, and thus $\Gamma^T\vdash_U (B[t/x])^T$, as desired.\\
\end{itemize}
This completes the induction.
\end{prf}
From (the proof of) this theorem, we can deduce a number of results. Our first result says that under the assumption of \cref{prov_in_U_and_U^T}, our construction does not build inconsistent theories out of consistent ones.
\begin{cor}\label{Sigma_and_U^T}
If $U$ and $T$ are theories such that $\vdash_U A$ implies $\vdash_U\Box_T A$ for all $\mathcal{L}$-formulae $A$, then the theories $U$ and $U^T$ prove the same $\Sigma_1$-formulae. In particular, $U^T$ is consistent if and only if $U$ is consistent.
\end{cor}
\begin{prf}
The first statement follows immediately from \cref{Sigma_and_T} and \cref{prov_in_U_and_U^T}. The second statement now follows since $\bot\in\Sigma_1$.
\end{prf}
Our next corollary tells us that the $T$-translation respects equivalence over $\ea$.
\begin{cor}\label{proofs_under_T}
Let $T$ be a theory and let $A$ and $B$ be $\mathcal{L}$-formulae. If $A\vdash_\ea B$, then $A^T\vdash_\ea B^T$.
\end{cor}
\begin{prf}
If $A\vdash_\ea B$, then also $A\vdash_{\ea^T} B$. By applying \cref{prov_in_U_and_U^T} with $U\equiv\ea$, we find that $A^T\vdash_\ea B^T$.
\end{prf}
The following result is the formalized counterpart of \cref{prov_in_U_and_U^T}.
\begin{cor}\label{formalized_prov_in_U_and_U^T}
Let $U$, $V$ and $T$ be theories such that $\vdash_V \bew_U(x)\to\Box_U\bew_T(x)$. Then $\vdash_V \Box_{U^T}A\leftrightarrow \Box_U A^T$ for all $\mathcal{L}$-formulae $A$.
\end{cor}
\begin{rem}
By \cref{inclusion_and_Boxrule}(ii), the requirement on $U$, $V$ and $T$ is satisfied when $U\leq T$. \ruitje
\end{rem}
\begin{prf}[Proof of \cref{formalized_prov_in_U_and_U^T}]
The `$\leftarrow$'-direction is immediate as it follows from the definition of $U^T$, and it does not need the requirement on $U$, $V$ and $T$. Concretely, we have 
\[\vdash_\ea \text{Form}(x)\wedge\mathsf{Bew}_U(x^T)\to\ax_{U^T}(x)\to \mathsf{Bew}_{U^T}(x).\] From this, the desired result follows.

For the `$\rightarrow$'-direction, we formalize the proof of the left-to-right direction of \cref{prov_in_U_and_U^T} inside $V$. We need that the statements of \cref{basic_facts_BoxT}, \cref{completeness_for_translations} and \cref{technical_stuff} are verifiable in $\ea$, hence in $V$. If we restrict the result to the case where $\Gamma$ is empty, we get $\vdash_V \mathsf{Bew}_{U^T}(x)\to\mathsf{Bew}_U(x^T)$, from which the desired result will follow.
\end{prf}
Finally, we discuss the relationship between $\ha$ and the $T$-translation.
\begin{cor}\label{HA_and_Ttranslation}
Let $U$, $V$ and $T$ be theories.
\begin{itemize}
\item[\textup{(}i\textup{)}]	If $\ha\subseteq U$, and $\vdash_U A$ implies $\vdash_U\Box_T A$ for all $\mathcal{L}$-formulae $A$, then $\ha\subseteq U^T$.
\item[\textup{(}ii\textup{)}]	If $\ha\leq U$ and $\vdash_V \bew_U(x)\to\Box_U\bew_T(x)$, then we have $\vdash_V \Box_\ha A\to \Box_{U^T} A$ for all $\mathcal{L}$-formulae $A$.
\end{itemize}
\end{cor}
\begin{prf}
In the proof of \cref{prov_in_U_and_U^T}, we have shown the following: if $A$ is the induction axiom for a certain formula $B$ and $U$ proves the induction axiom for $B^T$, then $\vdash_U A^T$. If $\ha\subseteq U$, then $U$ proves all induction axioms, so $U$ also proves $A^T$ for all induction axioms $A$. We can conclude that $\ha\subseteq U^T$.

Statement (ii) follows by formalizing this argument in $V$.
\end{prf}

Next, we isolate a special class of $\mathcal{L}$-formulae that behaves well with respect to the $T$-translation.
\begin{defi}\index[symb]{aa@$\mathcal{A}$}
The set $\mathcal{A}$ is the smallest set of $\mathcal{L}$-formulae such that
\begin{itemize}
\item[(i)]	$\mathcal{A}$ contains all atomic $\mathcal{L}$-formulae;
\item[(ii)]	$\mathcal{A}$ is closed under conjunction, disjunction, and both existential and universal quantification;
\item[(iii)]	if $S\in\Sigma_1$ and $A\in\mathcal{A}$, then $S\to A\in\mathcal{A}$. \ruitje
\end{itemize}
\end{defi}
\begin{lem}\label{preservation}
Let $T$ be a theory and let $A\in\mathcal{A}$. Then $\vdash_\ea A^T\to A$.
\end{lem}
\begin{prf}
We proceed by induction on the complexity of $A$. Only clause (iii) in the definition of $\mathcal{A}$ is nontrivial. Suppose that $A$ is $S\to B$, where $S\in\Sigma_1$, and that we already know the result for $B$. Then $\vdash_\ea S\leftrightarrow S^T$ and $\vdash_\ea B^T\to B$, so
\[\vdash_\ea (S\to B)^T\to (S^T\to B^T) \to (S\to B),\]
as desired.
\end{prf}
At the beginning of this section, we promised to construct theories that prove their own completeness. We now make this precise.

\begin{defi}\label{def_CP}
\index[term]{Completeness principle}\index[symb]{acp@CP}\index[term]{Strong L\"ob principle}\index[symb]{aslp@SLP}
Let $P(x)$ be a provability predicate for a theory $T$. Again, if $A$ is an $\mathcal{L}$-sentence, we write $\Box A$ for $P(\lg A\rg)$.
\begin{itemize}
\item[(i)]	The \emph{completeness principle} $\text{CP}_P$ is the axiom scheme $A\to\Box A$, where $A$ is an $\mathcal{L}$-sentence.
\item[(ii)]	The \emph{strong L\"ob principle} $\text{SLP}_P$ is the axiom scheme $(\Box A\to A)\to A$, where $A$ is an $\mathcal{L}$-sentence.
\end{itemize}
We will also write $\text{CP}_\Box$ instead of $\text{CP}_P$. This is actually a slight abuse of notation, because $\Box$ is merely an abbreviation and $\text{CP}_\Box$ really depends on $P(x)$. We write $\text{CP}_T$ for $\text{CP}_{\Box_T} = \text{CP}_{\bew_T}$. Similar conventions holds for SLP. \ruitje
\end{defi}
It turns out that the two schemes introduced above are two guises of the same principle.
\begin{lem}
Let $P(x)$ be a provability predicate for a theory $T$. Then the $\textup{CP}_P$ and $\textup{SLP}_P$ are interderivable over $\ea$.
\end{lem}
\begin{prf}
Define $\Box$ as above, and let $A$ be an $\mathcal{L}$-sentence. First, we show that $\vdash_{\ea+\text{CP}_\Box}\text{SLP}_\Box$. By \cref{properties_pp}(iii), we have
\[\vdash_{\ea+\text{CP}_\Box} (\Box A\to A)\to\Box(\Box A\to A)\to \Box A,\]
from which $\vdash_{\ea+\text{CP}_\Box} (\Box A\to A)\to A$ follows.

Now we show that $\vdash_{\ea+\text{SLP}_\Box} \text{CP}_\Box$. Clearly, we have $\vdash_\ea \Box(A\wedge \Box A)\to \Box A$, so
\begin{align*}
\vdash_{\ea+\text{SLP}_\Box} A &\to (\Box(A\wedge\Box A)\to A\wedge\Box A)\\
&\to A\wedge\Box A\\
&\to \Box A,
\end{align*}
as desired.
\end{prf}
Finally, here is the result we promised.
\begin{lem}\label{U^ast_proves_CP_U^ast}
For all theories $U$ and $T$, we have $\vdash_{U^T}\textup{CP}_{T^\ast}$. In particular, $\vdash_{U^\ast} \textup{CP}_{U^\ast}$.
\end{lem}
\begin{prf}
Let $A$ be an $\mathcal{L}$-sentence. Since $\Box_TA^T\in\Sigma_1$, we have that $\vdash_\ea \Box_TA^T\to(\Box_TA^T)^T$. So we get $\vdash_\ea A^T\to\Box_TA^T\to (\Box_TA^T)^T$, and also $\vdash_\ea \Box_T(A^T\to (\Box_TA^T)^T)$. So we find $\vdash_\ea (A\to\Box_TA^T)^T$. This means that $\vdash_U (A\to\Box_TA^T)^T$ as well, so we find that $\vdash_{U^T} A\to\Box_TA^T\to\Box_{T^\ast} A$, as desired. The second statement follows by taking $U\equiv T$.
\end{prf}
\begin{rem}
We remark that the proof of \ref{U^ast_proves_CP_U^ast} also goes through if we replace the first line with `Let $A$ be an $\mathcal{L}$-formula.' We will not need this greater generality. \ruitje
\end{rem}

\subsection{Fast and Slow Provability}\label{sec_fast_slow}

In this section, we introduce two nonstandard notions of provability. The first of these is \emph{fast provability}, which can be seen as iterated provability. The second is \emph{slow provability}, a notion of provability that puts a certain size restriction on the axioms that may be used in a proof. For developing the theory of fast provability, the following technique, that is also used in \cite{Paula}, will prove useful.

\begin{lem}[Reflexive induction]\index[term]{Reflexive induction}
Let $U\subseteq T$ be theories. Suppose $A(x)$ is an $\mathcal{L}$-formula in one free variable such that $\vdash_U A[0/x]$ and $\vdash_U \Box_T A\to A[\mathsf{S}x/x]$. Then $\vdash_U A$.
\end{lem}
\begin{prf}
It is provable in intuitionistic predicate logic that $\forall x\s A\to A$. So from our assumptions, it follows that $\vdash_U \Box_T\forall x\s A \to \Box_T A\to A[\mathsf{S}x/x]$. Since we also know that $\vdash_U A[0/x]$, we get $\vdash_U \Box_T\forall x\s A\to\forall x\s A$. Using L\"ob's Principle, we can conclude that $\vdash_U \forall x\s A$, so $\vdash_U A$.
\end{prf}

\begin{defi}\index[term]{Iterated provability}\index[symb]{aboxtf@$\Box^f_T$}\index[symb]{aboxtu@$\Box^{u+1}_T$}\index[symb]{aibewt@$\ibew_T$}\index[symb]{abewtf@$\bew^f_T$}
Let $T$ be a theory.
\begin{itemize}
\item[(i)]	We define $\ibew_T(u,x)$ as a formula satisfying
\[\vdash_\ea \ibew_T(u,x)\leftrightarrow ((u=0\wedge\bew_T(x))\vee\exists v\s (u=\mathsf{S}v\wedge \Box_T \ibew_T(v,x)))\]
as provided by the Diagonalization Lemma.
\item[(ii)]	For an $\mathcal{L}$-formula $A=A(x_1,\ldots,x_n)$, we write $\Box_T^{u+1}A$ for $\ibew_T(u,\lg A(\tilde{x}_1,\ldots,\tilde{x}_n)\rg)$
\item[(iii)]	We write $\bew^f_T(x)$ for $\exists u\s \ibew_T(u,x)$. Furthermore, for an $\mathcal{L}$-formula $A$, we write $\Box^f_TA$ for $\exists u\s \Box^{u+1}_TA$. \ruitje
\end{itemize}
\end{defi}
\begin{rem}\index[term]{Fast provability}
As we shall see shortly, $\bew^f_T$ is a provability predicate for $T$. This notion of provability is called \emph{fast provability} and was introduced by Parikh in \cite{parikh}. In this paper, fast provability is introduced in a different way, namely by closing the set of theorems of $T$ under Parikh's rule `from $\vdash\Box_TA$, infer $\vdash A$', where $A$ is an $\mathcal{L}$-sentence. This yields, verifiably in $\ha$, the same notion of provability we defined here. If $T$ is $\Sigma_1$-sound, then Parikh's rule does not lead to any new theorems, so the notions of ordinary provability and fast provability coincide. However, the use of Parikh's rule can lead to much shorter proofs, which explains the name `fast provability'. Later in this section, we will show that, if $T$ is consistent, it is never verifiable in $\ea$ that fast provability coincides with ordinary provability. \ruitje
\end{rem}

We notice that $\ibew_T$ is equivalent, over $\ea$, to a $\Sigma_1$-formula. Informally, $\ibew_T(u,x)$ can be thought of as the formula $\bew_T(\cdots(\bew_T(x))\cdots)$, where the $\bew_T$ occurs $u+1$ times, so we can see $\ibew_T$ as representing `iterated provability'. Notice that we write `$u+1$' in the superscript of $\Box_T$, to indicate that the $\Box_T$ `occurs' $u+1$ times. We prove a number of technical facts about $\Box^{u+1}_T$ and $\Box^f_T$.
\begin{lem}\label{facts_iterated_box}
Let $T$ be a theory and let $A$ and $B$ be $\mathcal{L}$-formulae. Then we have:
\begin{itemize}
\item[\textup{(}i\textup{)}]	$\vdash_\ea \Box_T\Box^{u+1}_TA\leftrightarrow\Box^{\mathsf{S}u+1}_TA\leftrightarrow\Box_T^{u+1}\Box_TA$;
\item[\textup{(}ii\textup{)}]	$\vdash_\ea u\leq v\to (\Box^{u+1}_TA\to\Box^{v+1}_TA)$;
\item[\textup{(}iii\textup{)}]	$\vdash_\ea \Box^{u+1}_T(A\to B)\to (\Box^{u+1}_TA\to\Box^{u+1}_TB)$,
\item[\textup{(}iv\textup{)}]	$\vdash_\ea \Box_TA\to\Box^f_TA$;
\item[\textup{(}v\textup{)}]	$\bew^f_T$ is a provability predicate for $T$;
\item[\textup{(}vi\textup{)}]	if $T$ is $\Sigma_1$-sound, then $\mathbb{N}\models \Box^f_TA$ if and only if $\mathbb{N}\models \Box_T A$, if and only if $\vdash_T A$;
\item[\textup{(}vii\textup{)}]	$\vdash_\ea \Box^f_T\Box_TA\leftrightarrow \Box^f_TA$.
\item[\textup{(}viii\textup{)}]	if $T$ is consistent, then $\nvdash_\ea \Box^f_T\bot\to\Box_T\bot$.
\end{itemize}
\end{lem}
\begin{prf}
(i)	From the definition of $\ibew_T$, it follows that $\vdash_\ea \Box_T\ibew_T(u,x)\leftrightarrow\ibew_T(\mathsf{S}u,x)$, so the first equivalence is immediate. For the second equivalence, we proceed by reflexive induction. First of all, we have $\vdash_\ea \Box^{\mathsf{S}0+1}_TA\leftrightarrow \Box_T\Box^{0+1}_TA \leftrightarrow \Box_T\Box_TA \leftrightarrow \Box^{0+1}_T\Box_TA$. Furthermore,
\begin{align*}
\vdash_\ea \Box_T(\Box^{\mathsf{S}u+1}_TA&\leftrightarrow\Box_T^{u+1}\Box_TA) \to\\
[\Box^{\mathsf{S}\mathsf{S}u+1}_TA &\leftrightarrow \Box_T\Box^{\mathsf{S}u+1}_TA\\
&\leftrightarrow \Box_T\Box^{u+1}_T\Box_TA\\
&\leftrightarrow \Box^{\mathsf{S}u+1}_T\Box_TA],
\end{align*}
which completes the proof.

(ii) Since $\Box^{u+1}_TA$ is a $\Sigma_1$-formula, we have $\vdash_\ea \Box^{u+1}_TA\to
\Box_T\Box^{u+1}_TA\to\Box^{\mathsf{S}u+1}_TA$. Now the claim follows by induction on $v$ inside $\ea$.

(iii) We proceed by reflexive induction. First of all, we have \[\vdash_\ea \Box^{0+1}_T(A\to B) \leftrightarrow \Box_T(A\to B) \to (\Box_TA\to\Box_TB) \leftrightarrow (\Box^{0+1}_TA\to \Box^{0+1}_TB).\]
Furthermore,
\begin{align*}
\vdash_\ea \Box_T(\Box^{u+1}_T(A\to B)&\to(\Box^{u+1}_TA\to\Box^{u+1}_TB))\to\\
[\Box^{\mathsf{S}u+1}_T(A\to B) &\leftrightarrow \Box_T\Box^{u+1}_T (A\to B)\\
&\to \Box_T(\Box^{u+1}_TA\to \Box^{u+1}_TB)\\
&\to (\Box_T\Box^{u+1}_TA\to\Box_T\Box^{u+1}_TB)\\
&\leftrightarrow (\Box^{\mathsf{S}u+1}_TA\to \Box^{\mathsf{S}u+1}_TB)],
\end{align*}
which completes the proof.

(iv) This is immediate as $\vdash_\ea \Box_TA\leftrightarrow\Box^{0+1}_TA$.

(v) This follows easily from (ii), (iii) and (iv).

(vi) The second equivalence was already asserted in \cref{basic_facts_BoxT}(i). So we prove the first equivalence. The right-to-left direction follows from (iv). For the converse, suppose that $\mathbb{N}\models\Box^{n+1} A$ for a certain $n\in\mathbb{N}$. If $n=0$, then we are done. So suppose that $n=m+1$ for a certain $m\geq 0$. Then $\mathbb{N}\models\Box_T\Box_T^{m+1}A$, so $\vdash_T\Box^{m+1}_TA$. Since $T$ is $\Sigma_1$-sound, we see that $\mathbb{N}\models \Box^{m+1}_TA$. By repeating this argument, we find $\mathbb{N}\models \Box_TA$, as desired.

(vii) This follows from $\vdash_\ea \Box^{u+1}_T\Box_T A\to \Box_T^{\mathsf{S}u+1}A$ and $\vdash_\ea \Box_T^{u+1}A\to \Box_T^{\mathsf{S}u+1} \to \Box^{u+1}_T\Box_TA$.

(viii) Suppose that $\vdash_\ea\Box^f_T\bot\to\Box_T\bot$. Then by (iv) and (vii), we have
\[\vdash_\ea \Box_T\Box_T\bot \to \Box^f_T\Box_T\bot \to\Box^f_T\bot \to\Box_T\bot,\]
so by L\"ob's Principle, we get $\vdash_\ea\Box_T\bot$. Since $\ea$ is $\Sigma_1$-sound, we conclude that $\vdash_T\bot$.
\end{prf}
\begin{rem}
It seems that in item (ii) above, we really need the presence of induction over $\Sigma_1$-formulae. The formula $\Box_T^{u+1}A\to\Box_T^{\mathsf{S}u+1}$ even provable in an intuitionistic version of Elementary Arithmetic. This means that passing from a witness of $\Box_T^{u+1}A$ to a witness of $\Box^{\mathsf{S}u+1}_TA$ is quite manageable, since this process is bounded by an elementary function. When producing a witness of $\Box_T^{v+1}A$ from a witness of $\Box_T^{u+1}A$, however, we need to iterate this process $v-u$ times, which means that the bound becomes a lot larger: possibly too large for weaker theories to handle. \ruitje
\end{rem}
Now we prove the analogue of \cref{formalized_prov_in_U_and_U^T} for fast provability.
\begin{lem}\label{formalized_prov_in_U^T_and_U_fast}
Suppose $U\leq T$ are theories and let $A$ be an $\mathcal{L}$-formula. Then we have $\vdash_\ea \Box_{U^T}^{u+1} A \leftrightarrow \Box_U^{u+1} A^T$. In particular, $\vdash_\ea \Box^f_{U^T}A\leftrightarrow\Box^f_U A^T$
\end{lem}
\begin{prf}
We proceed by reflexive induction. First of all, by \cref{formalized_prov_in_U_and_U^T}, we have
\[\vdash_\ea \Box_{U^T}^{0+1}A \leftrightarrow \Box_{U^T}A\leftrightarrow \Box_UA^T \leftrightarrow\Box_U^{0+1}A^T.\]
Furthermore, by \cref{Sigma_and_T}, \cref{formalized_prov_in_U_and_U^T} and \cref{facts_iterated_box}(vii), we have
\begin{align*}
\vdash_\ea \Box_U(\Box_{U^T}^{u+1} A &\leftrightarrow \Box_U^{u+1} A^T)\to\\
[\Box_{U^T}^{\mathsf{S}u+1} A &\leftrightarrow \Box_{U^T}\Box_{U^T}^{u+1}A\\
&\leftrightarrow \Box_U(\Box_{U^T}^{u+1}A)^T\\
&\leftrightarrow \Box_U\Box_{U^T}^{u+1}A\\
&\leftrightarrow \Box_U\Box_U^{u+1} A^T\\
&\leftrightarrow \Box_U^{\mathsf{S}u+1} A^T].
\end{align*}
This completes the proof.
\end{prf}

Now we turn to \emph{slow provability}.\index[term]{Slow provability} We will not give as many details as we did for fast provability, but instead we will refer to the paper \cite{PaulaFedor}. There are two reasons for this. First of all, developing the theory of slow provability is rather involved, so reasons of space do not permit us to provide all the details. The second reason involves our intended usage of fast and slow provability. In \cref{chap_applications}, we will obtain results about the provability logic of fast provability. In order to understand and appreciate these results, it is important to know what fast provability is, exactly. Slow provability, on the other hand, will only be used as a tool to obtain results that themselves do not mention slow provability. In order to understand these results, it is not necessary to know all the details about slow provability.

In the paper \cite{PaulaFedor}, the authors define a certain `fast-growing' total recursive function $F\colon \mathbb{N}\to\mathbb{N}$.\index[term]{Fast growing@Fast-growing function} There exists a $\Sigma_1$-formula $\varphi_F(x,y)$ representing $F$ in $\ha$. This means that the definition of $F$ as a recursive function is verifiable in $\ha$, and we have
\[\vdash_\ha \forall y\s (\varphi_F(n,y)\leftrightarrow y=F(n)) \quad\mbox{for all }n\in\mathbb{N}.\]
The $\Sigma_1$-formula $F(x)\!\!\downarrow$, which we read as `$F(x)$ is defined', is shorthand for $\exists y\s \varphi_F(x,y)$. We clearly have that $\vdash_\ha F(n)\!\!\downarrow$ for all $n\in\mathbb{N}$. However, the fast-growing function $F$ is constructed is such a way that $F$ is not provably total. That is, we do \emph{not} have $\vdash_\ha F(x)\!\downarrow$. Now we are ready to define slow provability.
\begin{defi}\index[term]{Heyting Arithmetic!Slow Heyting Arithmetic}\index[symb]{asha@$\sha$}
The theory \emph{slow Heyting Arithmetic}, denoted $\sha$, is given by the axiom formula
\[\ax_\sha(x):\leftrightarrow \ax_\ea(x)\vee(\ax_\ha(x)\wedge F(x)\!\downarrow). \wruitje\]
\end{defi}
Intuitively, we demand that the axioms we use must not be `too large': they must not be so large that they are beyond the domain of $F$. Since $F$ is in fact total, we see that $\mathbb{N}\models \ax_\sha(x)\leftrightarrow\ax_\ha(x)$, which means that $\ha=\sha$. We also clearly have that $\vdash_\ea \ax_\sha(x)\to\ax_\ha(x)$, so $\sha\leq\ha$. However, as we shall show shortly, we do \emph{not} have $\ha\leq\sha$. So from the viewpoint of $\ha$, the requirement that the axioms must not be too large is a genuine one.

Even though the base theory used in the paper \cite{PaulaFedor} is the classical theory $\pa$, many results carry over to the present case. The most important of these is:
\begin{prop}\label{prov_in_hahat}
We have $\vdash_\ha \bew_\ha(x)\to\Box_\ha\bew_\sha(x)$, and in particular, we have $\vdash_\ha \Box_{\ha^\sha} A\leftrightarrow \Box_\ha A^\sha$ for all $\mathcal{L}$-formulae $A$.
\end{prop}
\begin{prf}
The first statement is proven as in \cite{PaulaFedor}, Corollary 15, taking $\mathsf{S}_n$ to be the theory axiomatized by the axioms of $\ha$ having G\"odel number at most $n$. The second statement follows from \cref{formalized_prov_in_U_and_U^T} with $U\equiv V\equiv\ha$ and $T\equiv\sha$.
\end{prf}
The converse of this result, which is valid for the classical case, does \emph{not} carry over to the current setting, because the authors of \cite{PaulaFedor} use a model theoretic argument to derive this result. However, we will only need a very weak version of this converse, which we can `steal' from the classical case.
\begin{prop}\label{stealing_from_PA}
\begin{itemize}
\item[\textup{(}i\textup{)}]	For all $\Sigma_1$-sentences $S$, we have $\vdash_\ha \Box_\ha \Box_\sha S\to \Box_\ha S$.
\item[\textup{(}ii\textup{)}]	We have $\nvdash_\ha \Box_\ha\bot\to\Box_\sha\bot$. In particular, $\ha\nleq\sha$.
\end{itemize}
\end{prop}
\begin{prf}
(i) We define the analogue of slow provability for $\pa$, e.g.\@ by setting\index[symb]{asha1@$\mathsf{sPA}$}
\[\ax_\mathsf{sPA}(x):\leftrightarrow \ax_\ea(x)\vee(\ax_\pa(x)\wedge \exists y\geq x\s (F(y)\!\downarrow)).\]
Since $\vdash_\ea \ax_\ha(x)\to\ax_\pa(x)$, it is clear that $\ha\leq\pa$ and $\sha\leq\mathsf{sPA}$. We know from \cite{PaulaFedor}, Theorem 4, that $\vdash_\pa \Box_\pa\Box_\mathsf{sPA} S\to \Box_\pa S$. So we get
\[\vdash_\pa \Box_\ha\Box_\sha S \to \Box_\pa\Box_\mathsf{sPA} S\to\Box_\pa S\to \Box_\ha S,\]
where the final step holds since $\pa$ is, verifiably in $\ha$, $\Sigma_1$-conservative over $\ha$. We notice that $\Box_\ha\Box_\sha S \to \Box_\ha S$ is equivalent, over $\ha$, to a $\Pi_2$-sentence, that is, a sentence of the form $\forall x\s R(x)$, where $R\in\Sigma_1$. Since $\pa$ is $\Pi_2$-conservative over $\ha$, we also find that $\vdash_\ha \Box_\ha\Box_\sha S \to \Box_\ha S$, as desired.

(ii) Suppose that $\vdash_\ha \Box_\ha\bot\to\Box_\sha\bot$. Since $\bot\in\Sigma_1$, we have
\[\vdash_\ha \Box_\ha\Box_\ha\bot \to \Box_\ha\Box_\sha\bot \to \Box_\ha\bot,\]
so by L\"ob's Theorem, we get $\vdash_\ha \Box_\ha\bot$. But then $\ha$ is inconsistent, contradiction. For the second statement, we observe that $\ha\leq\sha$ would imply that $\Box_\ha\bot\to\Box_\sha\bot$ is provable in $\ea$, hence also in $\ha$.
\end{prf}

 \begin{rem}
 There is an alternative approach to slow provability suggested by Fedor Pakhomov in conversation to Albert Visser. In this approach we can prove the analogue of Proposition~\ref{stealing_from_PA}(i)
 without the detour over {\sf PA} and without the restriction to $\Sigma_1$-sentences. See \cite{absorption}. \ruitje
 \end{rem}

\subsection{Intuitionistic Modal Logic}\label{sec_IML}

In this section, we briefly review intuitionistic modal logic, abbreviated IML, and we define the system of IML that will be relevant to us.\index[term]{Intuitionistic modal logic}\index[term]{Intuitionistic modal logic!IML} The language $\mathcal{L}_\Box$ of IML has a countable set of propositional constants, the absurdity sign $\perp$, the usual binary connectives $\wedge$, $\vee$ and $\to$, and the unary sentential operator $\Box$.\index[term]{Language!of IML}\index[symb]{albox@$\mathcal{L}_\Box$}\index[symb]{abox@$\Box$!as a modal operator} We shall also use $\mathcal{L}_\Box$ to denote the set of all $\mathcal{L}_\Box$-sentences. As our proof system, we pick a Hilbert-style system that has two inference rules:
\begin{center}
\AxiomC{$A$}
\AxiomC{$A\to B$}
\RightLabel{$\to$E}\BinaryInfC{$B$}
\DisplayProof \hspace{1cm} and \hspace{1cm}
\AxiomC{$A$}
\RightLabel{Nec}\UnaryInfC{$\Box A$}
\DisplayProof.
\end{center}\index[symb]{anec@Nec}
\begin{defi}
\begin{itemize}
\item[(i)]	The set $\mathsf{iK}\subseteq\mathcal{L}_\Box$  is the smallest set that contains:\index[symb]{aik@$\mathsf{iK}$}
\begin{itemize}
\item[(a)]	all ($\mathcal{L}_\Box$-substitution instances of) tautologies of intuitionistic propositional logic;
\item[(b)]	all $\mathcal{L}_\Box$-sentences of the form $\Box(A\to B)\to (\Box A\to\Box B)$, where $A,B\in\mathcal{L}_\Box$,\vspace{.25\baselineskip}
\end{itemize}
and is closed under $\to$E and Nec.
\item[(ii)]	A \emph{theory for IML} will be a set $T$ that satisfies $\mathsf{iK}\subseteq T\subseteq\mathcal{L}_\Box$ and is closed under $\to$E and Nec. If $A\in\mathcal{L}_\Box$ and $\Gamma\subseteq\mathcal{L}_\Box$, we write $\Gamma\vdash_T A$ if there exists a finite subset $\Gamma_0\subseteq\Gamma$ such that $\bigwedge\Gamma_0\to A$ is in $T$.\index[symb]{avdasht@$\vdash_T$!for a theory for IML}\index[term]{Theory!for IML}
\item[(iii)]	The theory $\mathsf{iGL}$ is the smallest theory for IML that contains $\mathsf{iK}$ and all sentences of the form $\Box(\Box A\to A)\to\Box A$, where $A\in\mathcal{L}_\Box$.\index[symb]{aigl@$\mathsf{iGL}$}
\item[(iv)]	The theory $\iGLC$ is the smallest theory for IML that contains $\mathsf{iGL}$ and all sentences of the form $A\to \Box A$, where $A\in\mathcal{L}_\Box$.\index[symb]{aiglc@$\iGLC$} \ruitje
\end{itemize}
\end{defi}\
We now proceed to define the semantics of intuitionistic modal logic.
\begin{defi}
\begin{itemize}
\item[(i)]	Consider a triple $\langle W,\int,\r\rangle$, where $W$ is a nonempty set and $\int$ and $\r$ are binary relations on $W$. We say that this triple satisfies the \emph{model property} if $\int\circ \r$ is a subrelation of $\r$. That is, for all $w,v,u\in W$ we should have: if $w\int v\r u$, then $w\r u$.\index[term]{Model property}
\item[(ii)]	A \emph{frame for IML} is a triple $\langle W, \int,\r\rangle$, where $W$ is a nonempty set and $\int $ and $\r$ are binary relations on $W$, such that: $\langle W,\int\rangle$ is a poset and $\langle W,\int,\r\rangle$ satisfies the model property.\index[term]{Frame for IML}
\item[(iii)]	A \emph{model for IML} is a quadruple $\langle W,\int,\r, V\rangle$, where $\langle W,\int,\r\rangle$ is a frame for IML and $V$ is a relation (called the \emph{valuation}) between $W$ and the proposition letters from $\mathcal{L}_\Box$ satisfying:
\[w\int v\mbox{ and }wVp\mbox{ implies }vVp, \]
for all $w,v\in W$ and proposition letters $p$.\index[term]{Model for IML}
\item[(iv)]	Let $M=\langle W,\int,\r,V\rangle$ be a model for IML, let $w\in W$ and let $A\in\mathcal{L}_\Box$. We define the forcing relation $M,w\Vdash A$ by recursion on $A$, as follows. For all $B,C\in \mathcal{L}_\Box$, we set:\index[term]{Forcing relation}\index[symb]{aforces@$\Vdash$}
\begin{itemize}
\item[(a)]	$M,w\Vdash p$ iff $wVp$ for all proposition letters $p$;
\item[(b)]	$M,w\Vdash B\wedge C$ iff $M,w\Vdash B$ and $M,w\Vdash C$;
\item[(c)]	$M,w\Vdash B\vee C$ iff $M,w\Vdash B$ or $M,w\Vdash C$;
\item[(d)]	$M,w\Vdash B\to C$ iff for all $v\in W$ such that $w\int v$ and $M,v\Vdash B$, we have $M,v\Vdash C$;
\item[(e)]	$M,w\Vdash \Box B$ iff for all $v\in W$ such that $w\r v$, we have $M,v\Vdash B$.\vspace{.5\baselineskip}
\end{itemize}
If $M$ is understood, we just write $w\Vdash A$ instead of $M,w\Vdash A$. We write $M\Vdash A$ if $M,w\Vdash A$ for all $w\in W$, in which case we say that $A$ is \emph{valid} on $M$. Given a frame $\langle W,\int,\r\rangle$ for IML, we say that $A\in\mathcal{L}_\Box$ is valid on this frame iff for all models $M=\langle W,\int,\r,V\rangle$ for IML, we have that $A$ is valid on $M$.\index[term]{Valid} \ruitje
\end{itemize}
\end{defi}
Usually, one writes `$R$' for the modal relation we call `$\r$' here. Our notation has certain advantages that will become apparent in the next section. We impose the model property on our frames because we want the following result:
\begin{prop}[Preservation of Knowledge]\index[term]{Preservation of Knowledge}
Let $M=\langle W,\int,\r,V\rangle$ be a model for IML. If we have $w,v\in W$ and $A\in\mathcal{L}_\Box$ such that $w\Vdash A$ and $w\int v$, then $v\Vdash A$.
\end{prop}
\begin{prf}
We proceed by induction on the complexity of $A$. The base case and the induction steps for conjunction, disjunction and implication are trivial. So suppose that $A$ is $\Box B$ and that we have $w,v\in W$ such that $w\int v$ and $w\Vdash \Box B$. Consider any $u\in W$ such that $v\r u$. Then $w\int v\r u$, so since $\langle W,\int,\r\rangle$ has the model property, we get $w\r u$. Since $w\Vdash\Box B$, it follows that $u\Vdash B$. Since $u$ was arbitrary, we can conclude that $v\Vdash \Box B$, as desired.
\end{prf}
For our purposes, the relevant frame properties are the following.
\begin{defi}
Let $\langle W,\int,\r\rangle$ be a frame for IML.
\begin{itemize}\index[term]{Transitive}\index[term]{Semi-transitive}\index[term]{Realistic}\index[term]{Conversely well-founded}\index[term]{Irreflexive}
\item[(i)]	We say that this frame is \emph{irreflexive} if $\r$ is irreflexive, that is, if $\neg(w\r w)$ for all $w\in W$.
\item[(ii)]	We say that this frame is \emph{transitive} if $\r$ is transitive, that is, if $\r\circ\r$ is a subrelation of $\r$.
\item[(iii)]	We say that this frame is \emph{semi-transitive} if $\r\circ\r$ is a subrelation of $\r\circ\int$.
\item[(iv)]	We say that this frame is \emph{realistic} if $\r$ is a subrelation of $\int$.
\item[(v)]	We say that this frame is \emph{conversely well-founded} if $\r$ is conversely well-founded, that is, if every nonempty subset of $W$ has a maximal element w.r.t.\@ $\r$.
\end{itemize}
We say that a model for IML has one of the properties mentioned above if the underlying frame has it.\ruitje
\end{defi}
The terminology from (iii) is not standard and was suggested by R.\@ Iemhoff. The idea behind it is as follows. We can view $\r$ as an accessibility relation that is relative to the various worlds, while $\int$ represents the `real' accessibility between worlds. If, in a realistic frame, a world $w$ thinks that some world $v$ is accessible, then $v$ is also \emph{really} accessible from $w$. We observe that, due to the model property, a realistic frame is automatically transitive. Indeed, suppose that $\langle W,\int,\r\rangle$ is a realistic frame for IML and suppose we have $w,v,u\in W$ such that $w\r v\r u$. Then we also have $w\int v\r u$, so $w\r u$ follows, as desired.

Now we relate our frame properties to the axioms of $\iGLC$.

\begin{prop}\label{frame_properties}
Let $F=\langle W,\int,\r\rangle$ be a frame for intuitionistic modal logic.
\begin{itemize}
\item[\textup{(}i\textup{)}]	The sentence $\Box(\Box p\to p)\to \Box p$ is valid on $F$ if and only if $F$ is semi-transitive and conversely well-founded.
\item[\textup{(}ii\textup{)}]	The sentence $p\to\Box p$ is valid on $F$ if and only if $F$ is realistic.
\end{itemize}
In particular, all theorems of $\iGLC$ are valid on all realistic and conversely well-founded frames.
\end{prop}
\begin{prf}
(i)	This result is known from the literature. We refer the reader to the paper \cite{modalelogica}, Lemma 8.

(ii) First, suppose that $F$ is realistic. Let $V$ be a valuation on $F$, and suppose we have $w\in W$ such that $w\Vdash p$. If $v\in W$ is such that $w\r v$, then also $w\int v$, so by preservation of knowledge, we get $v\Vdash p$. We conclude that $w\Vdash\Box p$, and thus that $p\rightarrow\Box p$ is valid on $F$.

Conversely, suppose that $F$ is not realistic. Then there exist $w,v\in K$ such that $w\r v$, but also $w\not\preceq v$. We define a valuation $V$ on $F$ such that
\[xVp \quad\mbox{if and only if}\quad w\int x.\]
Then $wVp$, but since $w\r v$ and $\neg(vVp)$, we also have $w\nVdash\Box p$. We conclude that $w\nVdash p\rightarrow\Box p$ and thus that $p\rightarrow\Box p$ is not valid on $F$.

The final statement is easily proven by an induction on $\iGLC$-proofs.
\end{prf}
In order to get a completeness theorem, we need the following terminology.
\begin{defi}
Let $T$ be a theory for intuitionistic modal logic.
\begin{itemize}
\item[(i)]	A set $X\subseteq \mathcal{L}_\Box$ is called \emph{adequate} if it is closed under taking subsentences.\index[term]{Adequate set}\index[term]{XSaturated@$X$-saturated set}
\item[(ii)]	Suppose $X\subseteq \mathcal{L}_\Box$ is adequate. A set $S\subseteq X$ is called \emph{$X$-saturated} if the following hold:\vspace{2pt}
\begin{itemize}
\item[(a)]	$S$ is consistent, that is, $S\nvdash_T\bot$;
\item[(b)]	if $A\in X$ and $S\vdash_T A$, then $A\in S$;
\item[(c)]	if $A\vee B\in S$, then $A\in S$ or $B\in S$. \ruitje
\end{itemize}
\end{itemize}
\end{defi}
Notice that the converse of item (b) also holds: if $A\in S$, then clearly $A\in X$ and $S\vdash_T A$. We will need the following result.
\begin{lem}[Extension Lemma]\label{extension_lemma}\index[term]{Extension Lemma}
Let $T$ be a theory for intuitionistic modal logic and let $X\subseteq\mathcal{L}_\Box$ be an adequate set. Suppose we have $R\subseteq X$ and $A\in\mathcal{L}_\Box$ such that $R\nvdash_T A$. Then there exists an $X$-saturated set $S\supseteq R$ such that $S\nvdash_T A$.
\end{lem}
\begin{prf}
We fix an enumeration $B_0, B_1, B_2, \ldots$ of the formulae in $X$ such that every element of $X$ occurs infinitely many times in the enumeration. We define the sequence $S_0\subseteq S_1\subseteq S_2\subseteq\ldots$ by recursion. First of all, we set $S_0=R$. Now suppose that $S_n$ has been defined. If $S_n\nvdash_T B_n$, then $S_{n+1}$ is just $S_n$. If $S_n\vdash_T B_n$, then
\[S_{n+1} = \begin{cases}
S_n\cup\{B_n\} &\mbox{if }B_n\mbox{ is not a disjunction};\\
S_n\cup\{B_n,C\} &\mbox{if }B_n\mbox{ is }C\vee D,\mbox{ and }S_n\cup\{C\}\nvdash_T A;\\
S_n\cup\{B_n,D\} &\mbox{if }B_n\mbox{ is }C\vee D,\mbox{ and }S_n\cup\{C\}\vdash_T A;\\
\end{cases}\]
We define $S$ as $\bigcup_{n\in\mathbb{N}} S_n$. Clearly, we have $S_n\subseteq X$ for all $n\in\mathbb{N}$, so $S\subseteq X$. 

Now we use induction on $n$ to prove that $S_n\nvdash_T A$ for all $n\in\mathbb{N}$. For $n=0$, this holds by assumption. Now suppose that $S_n\nvdash_T A$ for a certain $n\in\mathbb{N}$; we need to show that $S_{n+1}\nvdash_T A$. If $S_n\nvdash_T B_n$, then this holds trivially. So suppose that $S_n\vdash_T B_n$. Then we must have that $S_n\cup\{B_n\}\nvdash_T A$, so if $B_n$ is not a disjunction, then we are also done. So suppose that $B_n$ is $C\vee D$. If $S_n\cup\{C\}\nvdash_T A$, then we also have $S_{n+1} = S_n\cup\{B_n,C\}\nvdash_T A$, so we are done. Finally, suppose that $S_n\cup\{C\}\vdash_TA$. Then we cannot have $S_n\cup\{D\}\vdash_TA$. Indeed, if we have both $S_n\cup\{C\}\vdash_TA$ and $S_n\cup\{D\}\vdash_T A$, then also $S_n\cup\{C\vee D\}\vdash_TA$, which is not the case. So $S_n\cup\{D\}\nvdash_TA$, and it follows that $S_{n+1} = S_n\cup\{B_n,D\}\nvdash_TA$, as desired. This completes the induction.

It follows that $S\nvdash_TA$, and in particular, $S$ is consistent. We check that $S$ is $X$-saturated. Now suppose that $C\in X$ and $S\vdash_TC$. Then there must be an $n\in\mathbb{N}$ such that $S_n\vdash_TC$. Let $m\geq n$ be minimal such that $B_m$ is $C$. Then $S_m\vdash_TB_m$, so we get $B_m\in S_{m+1}\subseteq S$, that is $C\in S$. Finally, supppse that $C\vee D\in S$. Then there must be an $n\in\mathbb{N}$ such that $C\vee D\in S_n$. Let $m\geq n$ be minimal such that $B_m$ is $C\vee D$. Then $B_m\in S_n\subseteq S_m$, so we certainly have $S_m\vdash_T B_m$. It follows that $C\in S_{m+1}\subseteq S$ or $D\in S_{m+1}\subseteq S$. This concludes the proof.
\end{prf}
Using the Extension Lemma, we can prove a sound- and completeness theorem for $\iGLC$. This result also appears, in a stronger form, as Theorem 4.25 in \cite{ArdeshirMojtahedi}.
\begin{thm}\label{completeness_iGLC}
Let $A$ be an $\mathcal{L}_\Box$-sentence. Then $\vdash_\iGLC A$ if and only if $A$ is valid on all finite irreflexive realistic frames.
\end{thm}
\begin{prf}
It is well-known that any finite irreflexive transitive frame is conversely well-founded. So if $\vdash_\iGLC A$, then $A$ is indeed valid on all finite irreflexive realistic frames, by \cref{frame_properties}. Conversely, suppose that $\nvdash_\iGLC A$. Let $X_0$ be the set of subsentences of $A$, and let $X_1 = \{\Box B\mid B\in X_0\}$. Then $X:=X_0\cup X_1$ is an adequate set. We let $W$ be the set of all $X$-saturated sets. Clearly, $W$ is finite, and we have the subset relation $\subseteq$ on $W$. For $w,v\in W$, we write $w\r v$ if:\vspace{2pt}
\begin{itemize}
\item[(i)]	whenever $B\in\mathcal{L}_\Box$ and $\Box B\in w$, we have $B\in v$;
\item[(ii)]	there exists a $C\in\mathcal{L}_\Box$ such that $\Box C\not\in w$ and $\Box C\in v$.\\
\end{itemize}
For $w\in W$ and $p$ a proposition letter, we say that $wVp$ if and only if $p\in w$. We clearly have: if $wVp$ and $w\subseteq v$, then $vVp$. It is also not difficult to check that $\langle W,\subseteq,\r\rangle$ satisfies the model property. Finally, since $\nvdash_\iGLC A$, there exists a $w_0\in W$ such that $w_0\nvdash A$, by the Extension Lemma. In particular, $W$ is nonempty, so $M=\langle W,\subseteq,\r,V\rangle$ is a model for intuitionistic modal logic.

We claim that the frame $\langle W,\subseteq,\r\rangle$ is irreflexive and realistic. Irreflexivity is immediate from the definition. Now suppose we have $w,v\in W$ such that $w\r v$, and $B\in w$. If $B\in X_0$, then $\Box B\in X_1\subseteq X$ and $w\vdash_\iGLC \Box B$, so $\Box B\in w$. Since $w\r v$ we get $B\in v$. Now suppose that $B\in X_1$. Then $B$ is $\Box C$ for some $C\in X_0$. Since $w\r v$, we get $C\in v$. This means that $v\vdash_\iGLC B$, so $B\in v$. In both cases, we get $B\in v$, so we conclude that $w\subseteq v$, as desired.

Now we show that $w\Vdash B$ if and only if $B\in w$, for all $w\in W$ and $B\in X$. We proceed by induction on the complexity of $B$.
\begin{itemize}
\item[At]	For proposition letters, the result holds by the definition of $V$.
\item[$\wedge$]	Suppose that $B$ is $C\wedge D$ and that the result holds for $C$ and $D$. If $w\in W$, then $w\Vdash C\wedge D$ iff $w\Vdash C$ and $w\Vdash D$, iff $C\in w$ and $D\in w$. Now suppose that $C\in w$ and $D\in w$. Then $w\vdash_\iGLC C\wedge D$ and $C\wedge D\in X$, so $C\wedge D\in w$. Conversely, suppose that $C\wedge D\in w$. Then $w\vdash_\iGLC C,D$ and $C,D\in X$, so we get $C\in w$ and $D\in w$.
\item[$\vee$]	Suppose that $B$ is $C\vee D$ and that the result holds for $C$ and $D$. If $w\in W$, then $w\Vdash C\vee D$ iff $w\Vdash C$ or $w\Vdash D$, iff $C\in w$ or $D\in w$. Suppose that $C\in w$ or $D\in w$. Then in both cases, we have $w\vdash_\iGLC C\vee D$. Since $C\vee D\in X$, we get $C\vee D\in w$. Conversely, if $C\vee D\in w$, then $C\in w$ or $D\in w$ since $w$ is $X$-saturated.
\item[$\to$]	Suppose that $B$ is $C\to D$ and that the result holds for $C$ and $D$. If $w\in W$, then $w\Vdash C\to D$ iff for all $v\supseteq w$, we have that $v\Vdash C$ implies $v\Vdash D$. And this holds iff for all $v\supseteq w$, we have that $C\in v$ implies $D\in v$. Now suppose that $C\to D\in w$ and that we have $v\supseteq w$ such that $C\in v$. Then also $C\to D\in v$, so $v\vdash_\iGLC D$. Since $D\in X$, we get $D\in v$. Conversely, suppose that $C\to D\not\in w$. Since $C\to D\in X$, this means that $w\nvdash_\iGLC C\to D$, and hence $w\cup\{C\}\nvdash_\iGLC D$. Since $w\cup\{C\}\subseteq X$, we can use the Extension Lemma to find a $v\in W$ such that $w\cup\{C\}\subseteq v$ and $v\nvdash_\iGLC D$. Then $w\subseteq v$, $C\in v$, and $D\not\in v$, so it follows that $w\nVdash C\to D$.
\item[$\Box$]	Suppose that $B$ is $\Box C$ and that the result holds for $C$. If $w\in W$, then $w\Vdash \Box C$ iff for all $v\sqsupset w$, we have $v\Vdash C$. And this holds iff for all $v\sqsupset w$, we have $C\in v$. Now suppose that $\Box C\in w$ and that we have $v\sqsupset w$. Then by the definition of $\r$, we get $C\in v$. Conversely, suppose that $\Box C\not\in w$. Consider the set $R = \{D\in\mathcal{L}_\Box\mid\Box D\in w\}\cup\{\Box C\}\subseteq X$. Suppose that $R\vdash_\iGLC C$. Then $\{D\in\mathcal{L}_\Box\mid\Box D\in w\}\vdash_\iGLC \Box C\to C$, so we also get $\{\Box D\in\mathcal{L}_\Box\mid\Box D\in w\}\vdash_\iGLC \Box(\Box C\to C)$. In particular, $w\vdash_\iGLC \Box(\Box C\to C)$, which yields $w\vdash_\iGLC \Box C$. However, we also have $\Box C\in X$, so we get $\Box C\in w$, contradiction. So $R\nvdash_\iGLC C$. By the Extension Lemma, there exists a $v\in W$ such that $R\subseteq v$ and $v\nvdash_\iGLC C$. We have $\{D\in\mathcal{L}_\Box\mid\Box D\in w\}\subseteq v$, $\Box C\not\in w$ and $\Box C\in v$, so $w\r v$. Furthermore, we have $C\not\in v$, so $w\nVdash \Box C$.\\
\end{itemize}
This completes the induction. Since $w_0\nvdash_\iGLC A$, we have $A\not\in w_0$. Since $A\in X$, we can apply the above result to conclude that $w_0\nVdash A$. So $A$ is not valid on the finite irreflexive realistic frame $\langle W,\subseteq,\r\rangle$.
\end{prf}

\section{An Abstract Arithmetical Completeness Theorem}\label{chap_CT}

In this section, we prove a completeness theorem for certain kinds of provability logics. We prove the theorem in a rather abstract form, not yet mentioning any specific provability predicates. In \cref{sec_Solovay_function}, we introduce the general framework and define the required Solovay function along with the intended realization of the propositional letters of $\mathcal{L}_\Box$. \cref{sec_preservation} is of a rather technical nature and forms the heart of the proof. Here we show that the realization we defined commutes with the logical operators of $\mathcal{L}_\Box$. In \cref{sec_completeness_theorem}, we formulate the completeness theorem and use the preceding material to prove it.

\subsection{Definition of the Solovay Function}\label{sec_Solovay_function}

The general setting of this section is given by the following definition.
\begin{defi}\index[term]{Good pair}\index[symb]{atriangle@$\triangle$}
Let $T$ be a theory and let $P(x)$ and $Q(x)$ be $\Sigma_1$-formulae in one free variable. If $A$ is an $\mathcal{L}$-sentence, we write $\Box A$ for $P(\lg A\rg)$. We also write $\triangle A$ for $Q(\lg A\rg)$. We say that $(P,Q)$ is a \emph{good pair} for $T$ if the following conditions are satisfied:
\begin{itemize}
\item[(i)]	$P$ and $Q$ are provability predicates for $T$;
\item[(ii)]	if $\mathbb{N}\models \Box A$, then $\vdash_T A$, for all $\mathcal{L}$-sentences $A$;
\item[(iii)]	$\vdash_T\text{SLP}_\triangle$ (or equivalently, $\vdash_T\text{CP}_\triangle$);
\item[(iv)]	$\vdash_T \Box\triangle S\to\Box S$ for all $\Sigma_1$-sentences $S$. \ruitje
\end{itemize}
\end{defi}
We immediately observe that, if these clauses apply and $S$ is a $\Sigma_1$-sentence, then we also have $\vdash_T \triangle S\to\Box\triangle S\to \Box S$. We also notice that $\vdash_\ea A$ implies $\vdash_TA$, which implies $\vdash_\ea\Box A$ and $\vdash_\ea\triangle A$ for all $\mathcal{L}$-sentences $A$.
\begin{rem}
We remark that the definition of a good pair does not occur anywhere in the literature. This definition is extremely artificial and tailor made to obtain the result of this section. \ruitje
\end{rem}
In the remainder of this section, we suppose that a theory $T$ and a good pair $(P,Q)$ for $T$ are given. We also use $\Box$ and $\triangle$ as defined above.

Let $M_0=\langle W_0, \int_0, \r_0, V_0\rangle$ be a finite irreflexive realistic model for IML such that $W_0$ has a least element w.r.t.\@ $\int_0$. Let $r>0$ be the cardinality of $W_0$. We assume that $W_0=\{1,\ldots, r\}$ and that the node $r$ is the least element of $W_0$ w.r.t.\@ $\int_0$. Now we expand $M_0$ to a new model $M = \langle W,\int,\r,V\rangle$ for IML. Intuitively, we append a copy of $1+\omega^{\text{op}}$ (in the $\r$-order relation) to the node $r$. Formally, we do this as follows. We take $W = \mathbb{N} \supset W_0$. The relation $\int$ is defined by:
\begin{align*}
i\int j \quad\mbox{iff}\quad &1\leq i,j\leq r\mbox{ and }i\int_0 j,\\
&\mbox{or }i>r\mbox{ and }1\leq j\leq i,\\
&\mbox{or }i=0,
\end{align*}
for all $i,j\in\mathbb{N}$. The relation $\r$ is defined by:
\begin{align*}
i\r j \quad\mbox{iff}\quad &1\leq i,j\leq r\mbox{ and }i\r_0 j,\\
&\mbox{or }i>r\mbox{ and }1\leq j<i,\\
&\mbox{or }i=0\mbox{ and }j>0,
\end{align*}
for all $i,j\in\mathbb{N}$. Finally, $V$ is defined by:
\begin{align*}
iVp\quad\mbox{iff}\quad &1\leq i\leq r\mbox{ and }iV_0p,
\end{align*}
for all $i\in\mathbb{N}$ and proposition letters $p$.

We can prove that $M$ is again a realistic irreflexive model for IML; but of course $M$ is no longer finite. However, $M$ \emph{is} conversely well-founded, so $M$ still validates all theorems of $\mathsf{iGLC}$. Since $\int_0$ and $\r_0$ are finite relations, we can give $\Delta_0$-definitions of these relations inside $\ea$. Now we can formalize the definitions of $\int$ and $\r$ given above in order to obtain $\Delta_0$-definitions of $\int$ and $\r$ inside $\ea$. Then $\ea$ verifies the relevant properties of $M$: that $\langle \mathbb{N},\int\rangle$ is a poset, that $\r$ is irreflexive, that $\langle W,\int,\r\rangle$ has the model property, and that this frame is realistic. E.g.\@ by verification of the model property we mean that $\vdash_\ea x\int y\wedge y\r z\to x\r z$. Since $\int$ is defined by a $\Delta_0$-formula, we have: if $i\int j$, then $\vdash_\ea i\int j$, and if $i\npreceq j$, then $\vdash_\ea\neg(i\int j)$. A similar result holds for $\r$. Moreover, by \cref{Delta0_stuff}, we can safely make case distinctions like $x\int y\vee\neg(x\int y)$ inside $\ea$. Observe that these remarks also apply in $T$ instead of $\ea$.

For an $A\in\mathcal{L}_\Box$, we define the set $\lval A\rval$ as $\{i\in\mathbb{N}\mid i\Vdash A\}$. The model $M$ is constructed in such a way that the following result holds.\index[symb]{aat2@$\lval A\rval$}
\begin{lem}
If $A\in\mathcal{L}_\Box$, then $\lval A\rval$ is finite or $\lval A\rval =\mathbb{N}$.
\end{lem}
\begin{prf}
We have to show the following: if $i\in\lval A\rval$ for all $i>0$, then $0\in\lval A\rval$. We proceed by induction on the complexity of $A$. The atomic case clearly holds, and the steps for $\wedge$ and $\vee$ are trivial. Now suppose that $A$ is $B\to C$ and that the claim holds for $B$ and $C$. Suppose that $i\in\lval B\to C\rval$ for all $i>0$, and that $0\not\in\lval B\to C\rval$. Then we must have $0\in\lval B\rval$ and $0\not\in\lval C\rval$. By the induction hypothesis, $i\not\in C$ for some $i>0$. However, since $0\preceq i$, we also have $i\in\lval B\rval$, so $i\not\in\lval B\to C\rval$, contradiction. Finally, suppose that $A$ is $\Box B$ and that the claim holds for $B$. Suppose that $i\in\lval \Box B\rval$ for all $i>0$. We should show that $0\in\lval\Box B\rval$. By preservation of knowledge, it suffices to show that $j\in\lval B\rval$ for all $j\geq r$. But for such $j$, we have $j+1\in\lval \Box B\rval$ by assumption, and $j+1\sqsubset j$, so we indeed have $j\in\lval B\rval$.
\end{prf}
We now proceed to define the Solovay function. Our models are equipped with \emph{two} relations, as opposed to just one in the classical case, and we need to find some way to incorporate this into the Solovay function. We will use two separate provability predicates to take care of the relations $\int$ and $\r$. This is where our good pair comes in. Since $P(x)$ and $Q(x)$ are $\Sigma_1$-formulae, we can write $P(x)$ as $\exists y\s\pf_\Box(y,x)$ and $Q(x)$ as $\exists y\s\pf_\triangle(y,x)$, where $\pf_\Box$ and $\pf_\triangle$ are $\Delta_0$-formulae. \index[symb]{aprfbox@$\pf_\Box$}\index[symb]{aprfbox1@$\pf_\triangle$}

Let $\langle \cdot,\cdot\rangle\colon\mathbb{N}^2\to\mathbb{N}$ be a primitive recursive pairing function that can be formulated inside $\ea$ using a $\Delta_0$-formula. Let $p_0\colon\mathbb{N}\to\mathbb{N}$ be the elementary function that gives the projection onto the first coordinate. By replacing $\pf_\Box(y,x)$ with $\exists z\leq y\s (y=\langle x,z\rangle\wedge \pf_\Box(z,x))$, we may assume without loss of generality that
\begin{align}\label{unique_proofs}
\vdash_\ea \pf_\Box(y,x)\to x=p_0(y).
\end{align}
We do the same for $\pf_\triangle$.

In the sequel, we write $x\prec y$ for $x\int y\wedge\neg(x=y)$ and $x\req y$ for $x\r y\vee x=y$. We define the function $h\colon\mathbb{N}\to\mathbb{N}$ by $h(0)=0$ and\index[symb]{ah@$h$}\index[term]{Solovay function}
\[
h(k+1) = \begin{cases}
m\quad&\mbox{if } h(k)\sqsubset m\mbox{ and }\pf_\Box(k,\lg\exists x\s \neg(h(x)\sqsubseteq m)\rg);\\
n\quad&\mbox{if } h(k)\prec n\mbox{ and }\pf_\triangle(k,\lg\exists y\s \neg(h(y)\preceq n)\rg);\\
h(k)\quad&\mbox{if neither of these apply}.
\end{cases}\]
Here $x$ and $y$ are two (syntactically) distinct variables, so by our assumption \cref{unique_proofs} above, the first two clauses can never apply simultaneously. Using \cref{unique_proofs} again, we also see that $m$ as in the first clause, if it exists, is unique, and similarly for the second clause. Using the Diagonalization Lemma, we can give a $\Sigma_1$-definition of $h$ inside $\ea$. Then we can prove inside $\ea$ that $h$ is in fact a function. We also have $\vdash_\ea x\leq y\to h(x)\preceq h(y)$, which can be shown by induction on $y$ inside $\ea$.

Notice that it is in some sense `easier' to move along $\r$ than it is to move along $\int$. We have $\vdash_\ea \exists y\s\neg(h(y)\int m)\to \exists x\s\neg(h(x)\req m)$. Since $Q$ is a provability predicate for $T$, we also find that $\vdash_\ea \triangle(\exists y\s\neg(h(y)\int m))\to\triangle(\exists x\s\neg(h(x)\req m))$. We also observe that $\exists x\s\neg(h(x)\req m)$ is equivalent, over $\ea$, to a $\Sigma_1$-sentence. This means that we also have $\vdash_T \triangle(\exists x\s\neg(h(x)\req m))\to\Box(\exists x\s\neg(h(x)\req m))$. We conclude that
\begin{align}\label{moving_comparison}
\vdash_T \triangle(\exists y\s\neg(h(y)\int m))\to \Box(\exists x\s\neg(h(x)\req m)),
\end{align}
for any $m\in\mathbb{N}$. We will need this in the sequel. We also need the following observation: if $i\neq 0$ is a natural number, then
\begin{align}\label{trapped}
\vdash_\ea \neg(x\int i)\leftrightarrow\bigvee_{j\in U} x=j,
\end{align}
where $U = \{j\in\mathbb{N}\mid j\npreceq i\}$ is finite. In other words, if $\ea$ knows that $x\npreceq i$ for some standard $i\neq 0$, then $\ea$ knows that $x$ is some standard number as well. For $\req$, a similar remark applies.

We close this section with a definition.
\begin{defi}
For a sentence $A\in\mathcal{L}_\Box$, we define the $\mathcal{L}$-sentence $[A]$ as \index[symb]{aat1@$[A]$}
\[\begin{cases}
\bigvee_{i\in\lval A\rval} \exists x\s (h(x)=i) &\mbox{if }\lval A\rval\mbox{ is finite};\\
\top &\mbox{if }\lval A\rval=\mathbb{N}.
\end{cases}\wruitje\]
\end{defi}
We observe that $[A]$ is always (equivalent to) a $\Sigma_1$-sentence.

\subsection{Preservation of the Logical Structure}\label{sec_preservation}

In this rather technical section, we show that $[\cdot]$ commutes with all the logical operators figuring in $\mathcal{L}_\Box$. The proofs in this section will become increasingly difficult. We adopt all the notation introduced in the previous section.

\begin{lem}\label{vee_pres}
We have $\vdash_\ea [B\vee C]\leftrightarrow [B]\vee[C]$ for $B,C\in\mathcal{L}_\Box$.
\end{lem}
\begin{prf}
This is immediate from the definition of $[\cdot]$.
\end{prf}
\begin{lem}\label{wedge_pres}
We have $\vdash_\ea [B\wedge C]\leftrightarrow [B]\wedge[C]$ for $B,C\in\mathcal{L}_\Box$.
\end{lem}
\begin{prf}
If $\lval B\rval=\mathbb{N}$, then $\lval B\wedge C\rval=\lval C\rval$, so we have $\vdash_\ea [B\wedge C]\leftrightarrow [C] \leftrightarrow [B]\wedge[C]$. Similarly, the result follows if $\lval C\rval=\mathbb{N}$. So suppose that $\lval B\rval$ and $\lval C\rval$ are both finite; then $\lval B\wedge C\rval$ is finite as well.

The `$\rightarrow$'-statement is immediate in this case. For the other direction, we should show that $\vdash_\ea \exists x\s (h(x)=i)\wedge\exists y\s(h(y)=j)\to[B\wedge C]$ whenever $i\in\lval B\rval$ and $j\in\lval C\rval$. First of all, we notice that $\vdash_\ea \exists x\s (h(x)=i)\wedge\exists y\s (h(y)=j)\to i\int j\vee j\int i$. Indeed, reason inside $\ea$ and suppose we have $x$ and $y$ such that $h(x)=i$ and $h(y)=j$. Since $x\leq y\vee y\leq x$ and $x\leq y\to h(x)\int h(y)$, we can conclude that $i\int j\vee j\int i$, as desired.

Now, if $i$ and $j$ are incomparable w.r.t.\@ $\int$, then $\vdash_\ea\neg(i\int j\vee j\int i)$, so by the above we have $\vdash_\ea \neg(\exists x\s (h(x)=i)\wedge\exists y\s(h(y)=j))$, in which case the result is clear. If $i$ and $j$ are comparable w.r.t.\@ $\int$, then assume without loss of generality that $i\int j$. Then $j\in\lval B\wedge C\rval$, so 
\[\vdash_\ea \exists x\s (h(x)=i)\wedge\exists y\s(h(y)=j) \to \exists y\s(h(y)=j)\to[B\wedge C],\] as desired.
\end{prf}
\begin{lem}\label{to_pres}
We have $\vdash_T [B\to C]\leftrightarrow ([B]\to[C])$ for $B,C\in\mathcal{L}_\Box$.
\end{lem}
\begin{prf}
If $\lval B\to C\rval=\mathbb{N}$, then $\lval B\rval\subseteq\lval C\rval$, so $[B\to C]$ and $[B]\to[C]$ are both equivalent to $\top$, even over $\ea$. Now suppose that $\lval B\to C\rval$ is finite.

We first treat the $\leftarrow$-direction. Let $j_0, \ldots, j_{s-1}\neq 0$ be the $\preceq$-maximal elements $j$ of $\mathbb{N}$ such that $j\not\in\lval B\to C\rval$. Then for all $t<s$, we have $j_t\in\lval B\rval$ and $j_t\not\in\lval C\rval$. Using the fact that $\prec$ is also a conversely well-founded relation, we can show that for all $i\in\mathbb{N}$, we have $i\in\lval B\to C\rval$ if and only if $i\not\int j_t$ for all $t<s$.

Now we reason inside $\ea$. Suppose that $[B]\to[C]$ and $\triangle[B\to C]$. Since $[B\to C]\to\exists y\s \neg(h(y)\preceq j_t)$, we have $\triangle(\exists y\s\neg(h(y)\preceq j_t))$. Now let $k_t$ satisfy $\mathsf{Prf}_\triangle(k_t,\ulcorner \exists y\s\neg(h(y)\preceq j_t)\urcorner)$. We distinguish three cases (which is constructively acceptable).
\begin{enumerate}
\item	Suppose that $h(k_t)\prec j_t$ for some $t<s$. Then then by the definition of $h$, we get $h(k_t+1) = j_t$. But $j_t\in \lval B\rval$, so $[B]$ holds, so $[C]$ holds, and therefore $[B\to C]$ also holds.
\item	Suppose that $h(k_t)=j_t$ for some $t<s$. Then $[B\to C]$ again follows.
\item	Suppose that $\neg(h(k_t)\preceq j_t)$ for all $t<s$. Let $k=\max_{t<s} k_t$. Then we also know that $\neg(h(k)\preceq j_t)$ for all $t<s$. Indeed, suppose that $h(k)\preceq j_t$ for some $t$. Since $k_t\leq k$, we get $h(k_t)\preceq h(k)\preceq j_t$, so, since $\preceq$ is (provably) transitive, $h(k_t)\preceq j_t$, which we already excluded. So we indeed have $\neg(h(k)\preceq j_t)$ for all $t<s$. But then by \cref{trapped} applied to $j_0, \ldots, j_{s-1}$, we see that $\bigvee_{j\in U} h(k)=j$, where $U = \{j\in\mathbb{N}\mid j\not\int j_t\mbox{ for all }t<s\}$ is a finite set. We see (outside $\ea$) that $U= \lval B\to C\rval$, so (inside $\ea$ again) we get $[B\to C]$.\\
\end{enumerate}
We conclude that $\vdash_\ea ([B]\to[C])\to(\triangle[B\to C]\to[B\to C])$. Since $\vdash_T\text{SLP}_\triangle$, we may conclude that $\vdash_T ([B]\to[C])\to [B\to C]$.

The $\rightarrow$-direction is even provable in $\ea$. Notice that $\lval B\wedge(B\to C)\rval\subseteq \lval C\rval$, so by \cref{wedge_pres}, we have
\[\vdash_\ea ([B]\wedge [B\to C]) \to [B\wedge(B\to C)] \to [C],\]
so $\vdash_\ea [B\to C]\to([B]\to[C])$.
\end{prf}
Before we can show that $[\cdot]$ commutes with $\Box$, we first need some auxiliary results.
\begin{lem}\label{danger_of_leaving_home}
Suppose $i>0$ is a natural number. Then 
\[\vdash_T \exists x\s (h(x)=i)\to\Box(\exists y\s (i\prec h(y))).\]
\end{lem}
\begin{prf}
Before we start proving the displayed sentence inside $T$, we need to verify two auxiliary facts inside $\ea$. First of all, we claim that
\[\vdash_\ea(\neg(h(y)\preceq i)\wedge h(x)=i)\to i\prec h(y).\]
Reason inside $\ea$ and assume the antecedent. If $y<x$, then $h(y)\preceq h(x)=i$, quod non. So $x\leq y$, which means that $i=h(x)\preceq h(y)$. But $h(y)$ cannot be equal to $i$, so $i\prec h(y)$, as desired. Now we also have:
\begin{align}\label{case1}
\vdash_\ea(\exists y\s \neg(h(y)\preceq i)\wedge \exists x\s (h(x)=i))\to \exists y\s (i\prec h(y)).
\end{align}
Secondly, we claim that
\[\vdash_\ea(\neg(h(y)\sqsubseteq i)\wedge h(x)=i\wedge h(x-1)\r i) \to i\prec h(y).\]
Again, reason inside $\ea$ and assume the antecedent. Suppose that $y<x$. Then $y\leq x-1$, so $h(y)\preceq h(x-1) \sqsubset i$. Since our frame (provably) has the model property, we get $h(y)\sqsubset i$, contradiction. So $y\geq x$. But then $i=h(x)\preceq h(y)$ and $\neg(i=h(y))$, so $i\prec h(y)$, as desired. We also find: 
\begin{align}\label{case2}
\vdash_\ea\exists y\s (\neg(h(y)\sqsubseteq i))\wedge \exists x\s (h(x)=i\wedge h(x-1)\r i) \to \exists y\s (i\prec h(y)).
\end{align}
Now we start the main part of the proof. Reason inside $T$, and suppose that we have an $x$ such that $h(x)=i$. Since $h$ is (provably) a function, we can consider the \emph{least} $x$ such that $h(x)=i$. Then $x>0$, and $h(x-1)\prec i$. Again, we make a constructively acceptable case distinction.
\begin{enumerate}
\item	Suppose that $\neg (h(x-1)\sqsubset i)$. Then $\triangle(\exists y\s \neg(h(y)\preceq i))$ (otherwise, we wouldn't have moved up to $i$). Since $\exists x\s (h(x)=i)$ is a $\Sigma_1$-sentence, we also get $\triangle(\exists x\s (h(x)=i))$. Using \cref{case1} and the properties of $\triangle$, we can conclude that $\triangle\exists y\s(i\prec h(y))$. Since $\exists y\s (i\prec h(y))$ is a $\Sigma_1$-sentence, we also get $\Box(\exists y\s(i\prec h(y))$ by \cref{moving_comparison}, as desired.
\item	Suppose that $h(x-1)\sqsubset i$. Then, from the fact that we moved up to $i$, we can deduce that $\Box(\exists x\s\neg(h(x)\sqsubseteq i))$ or $\triangle(\exists y\s\neg(h(y)\preceq i))$. By \cref{moving_comparison}, we can conclude that $\Box(\exists y\s\neg(h(y)\sqsubseteq i))$ in both cases. At this point, we have $\exists x\s (h(x)=i\wedge h(x-1)\r i)$. Since this is a $\Sigma_1$-sentence, we also get $\Box(\exists x\s h(x)=i\wedge h(x-1)\r i))$. Using \cref{case2} and the properties of $\Box$, we again find $\Box\exists y\s(i\prec h(y))$, as desired.
\qedhere
\end{enumerate}
\end{prf}
\begin{lem}\label{running_away}
Let $i,j$ be natural numbers such that $i\prec j$ and $\neg(i\sqsubset j)$. Then 
\[\vdash_\ea \exists x\s (h(x)=i)\wedge\exists y\s (h(y)=j)\to \triangle(\exists z\s (j\prec h(z))).\]
\end{lem}
\begin{prf}
First of all, we notice that we also know that $i\prec j$ and $\neg(i\sqsubset j)$ inside $\ea$. Now reason inside $\ea$, and suppose that $\exists x\s (h(x)=i)$ and $\exists y\s (h(y)=j)$. Since $h$ is (provably) a function, we can consider the least $y$ such that $h(y)=j$. Then $y>0$, and $h(y-1)\prec j$. Consider an $x$ such that $h(x)=i$. Suppose that $y\leq x$. Then $j=h(y)\prec h(x)=i\prec j$, which is a contradiction since $\int$ is (provably) antisymmetric. So $x<y$, which also means $x\leq y-1$. Now we get $i=h(x)\preceq h(y-1)$.

If $h(y-1)\sqsubset j$, then $i\preceq h(y-1)\sqsubset j$, so $i\sqsubset j$. But we also have $\neg(i\sqsubset j)$, so we must have $\neg(h(y-1)\sqsubset j)$. Now we can use the exact same reasoning as in case 1 in the proof of \cref{danger_of_leaving_home} (with $j$ instead of $i$, and $y$ instead of $x$) to arrive at $\triangle\exists z\s (j\prec h(z))$, as desired. (Observe that we can perform this reasoning inside $\ea$ instead of $T$, since we do not need \cref{moving_comparison} here.)
\end{prf}
Now that we have proven these tedious lemmata, we can derive our crucial result.
\begin{lem}\label{box_pres}
We have $\vdash_T [\Box B]\leftrightarrow \Box[B]$ for all $B\in\mathcal{L}_\Box$.
\end{lem}
\begin{prf}
If $\lval B\rval=\mathbb{N}$, then $\lval \Box B\rval=\mathbb{N}$ as well, and we see that $[\Box B]$ and $\Box[B]$ are both equivalent to $\top$ over $\ea$. Now suppose that $\lval B\rval$ is finite.

We first treat the $\leftarrow$-direction, which can be shown even in $\ea$. Let $j_0, \ldots, j_{s-1}\neq 0$ be the $\sqsubset$-maximal elements $j$ of $\mathbb{N}$ such that $j\not\in \lval B\rval$. Notice that $j_t\in\lval\Box B\rval$ for all $t<s$. Suppose that we have $i\in\lval B\rval$ and $t<s$ such that $i\sqsubseteq j_t$. Since $M$ is realistic, we get $i\preceq j_t$, so by preservation of knowledge, $j_t\in\lval B\rval$, contradiction. So if $i\in\lval B\rval$, then $i\not\sqsubseteq j_t$. In particular, we have $\vdash_\ea[B]\to\exists x\s \neg(h(x)\sqsubseteq j_t)$ for all $t<s$. Using the transitivity of $\r$ and the fact that $\r$ is a conversely well-founded relation, we can also show: if $i\not\sqsubseteq j_t$ for all $t<s$, then $i\in \lval B\rval$.

Now we reason inside $\ea$ and suppose that $\Box[B]$. Then $\Box(\exists x\s \neg(h(x)\sqsubseteq j_t))$ also holds. Let $k_t$ satisfy $\mathsf{Prf}_\Box(k_t,\ulcorner\exists x\s \neg(h(x)\sqsubseteq j_t)\urcorner)$. We distinguish three cases.\vspace{.5\baselineskip}
\begin{enumerate}
\item 	Suppose $h(k_t)\sqsubset j_t$ for some $t<s$. Then by the definition of $h$, we have $h(k_t+1)=j_t$, and $[\Box B]$ follows.
\item	Suppose $h(k_t)=j_t$ for some $t<s$. Then $[\Box B]$ again follows.
\item	Suppose that $\neg(h(k_t)\sqsubseteq j_t)$ for all $t<s$. Let $k=\max_{t<s} k_t$. If $h(k)=j_t$ for some $t<s$, then $[\Box B]$ again follows. Suppose $h(k)\sqsubset j_t$ for some $t<s$. Since $h(k_t)\preceq h(k)\sqsubset j_t$ and our frame (provably) has the model property, we get $h(k_t)\sqsubset j_t$, which we already excluded. So we have $\neg(h(k)\sqsubseteq j_t)$ for all $t<s$. But then using the $\sqsubseteq$-analogue of \cref{trapped} for $j_0, \ldots, j_{s-1}$, we see that $\bigvee_{j\in U} h(k)=j$, where $U=\{j\in\mathbb{N}\mid j\not\sqsubseteq j_t\mbox{ for all }t<s\}$ is a finite set. We see (outside $\ea$) that $U = \lval B\rval\subseteq\lval \Box B\rval$, where the inclusion holds since $M$ is realistic. So (inside $\ea$ again), we get $[\Box B]$, as desired.\\
\end{enumerate}

Now we treat the $\rightarrow$-direction. Consider an $i\in\lval \Box B\rval$. Then $i>0$, since $\lval B\rval$ is finite. So by \cref{danger_of_leaving_home}, we have 
\begin{align}\label{the_first_step}
\vdash_T \exists x\s (h(x)=i)\to \Box(\exists y\s (i\prec h(y))).
\end{align}
Every nonzero node $k$ of $M$ has a finite $\prec$-rank, which is the greatest $n$ such that there exists a sequence $k=k_0\prec k_1\prec \cdots\prec k_n$. Let $a\in\mathbb{N}$ be the $\prec$-rank of $i$. For $b\in\mathbb{N}$, we define the finite set
\[U_b = \{j\in\mathbb{N}\mid i\prec j, i\not\r j\mbox{ and } \text{rank}(j)<b\}.\]
We know (inside $\ea$) that $i\prec h(y)$ implies that $h(y)$ is a standard number. Moreover, such a standard number must have rank smaller than $a$, so it is either in $\lval B\rval$ (if $i\r h(y)$) or in $U_a$ (if $i\not\r h(y)$). That is, we have
\[\vdash_\ea i\prec h(y) \to \bigvee_{j\in \lval B\rval} h(y)=j \vee\bigvee_{j\in U_a} h(y)=j.\]
From this, it follows that
\[\vdash_\ea \exists y\s(i\prec h(y))\to [B]\vee \bigvee_{j\in U_a}\exists y\s (h(y)=j).\]
So \cref{the_first_step} together with the properties of $\Box$ implies that
\begin{align}\label{leftovers}
\vdash_T \exists x\s (h(x)=i)\to \Box\left([B]\vee \bigvee_{j\in U_a}\exists y\s (h(y)=j)\right).
\end{align}
Suppose that $j\in U_b$ for a certain $b\geq 1$. By \cref{running_away}, we know that 
\begin{align}\label{the_second_step}
\vdash_\ea \exists x\s (h(x)=i)\wedge\exists y\s (h(y)=j)\to \triangle(\exists z\s (j\prec h(z))).
\end{align}
Furthermore, if $j\prec h(z)$, then we know (inside $\mathsf{HA}$) that $h(z)$ is some standard number. Moreover, such a standard number must have lower $\prec$-rank than $j$, so it is either in $\lval B\rval$ (if $i\r h(z)$) or in $U_{b-1}$ (if $i\not\r h(z)$). That is, we have 
\[\vdash_\ea j\prec h(z)\to \bigvee_{k\in\lval B\rval}h(z)=k\vee\bigvee_{k\in U_{b-1}} h(z)=k.\] 
From this, it follows that 
\[\vdash_\ea \exists z\s (j\prec h(z))\to[B]\vee\bigvee_{k\in U_{b-1}} \exists z\s (h(z)=k).\]
So using \cref{the_second_step} and the properties of $\triangle$, we get
\[\vdash_\ea \exists x\s (h(x)=i)\wedge\exists y\s(h(y)=j) \to \triangle\left([B]\vee\bigvee_{k\in U_{b-1}} \exists z\s(h(z)=k)\right).\]
This holds for all $j\in U_b$, so
\[\vdash_\ea \exists x\s (h(x)=i)\wedge\bigvee_{j\in U_b}(\exists y\s(h(y)=j)) \to \triangle\left([B]\vee\bigvee_{j\in U_{b-1}} \exists y\s(h(y)=j)\right).\]
(We changed some bound variables on the right hand side.)
Since $[B]$ is equivalent, over $\ea$, to a $\Sigma_1$-sentence, we also have
\[\vdash_\ea [B]\to \triangle[B] \to \triangle\left([B]\vee \bigvee_{j\in U_{b-1}} \exists y\s(h(y)=j)\right).\]
So we conclude that
\[\vdash_\ea \exists x\s (h(x)=i)\wedge\left([B]\vee\bigvee_{j\in U_b}(\exists y\s (h(y)=j))\right)\to \triangle\left([B]\vee \bigvee_{j\in U_{b-1}} \exists y\s(h(y)=j)\right).\]
Since $\exists x\s(h(x)=i)$ is equivalent, over $\ea$, to a $\Sigma_1$-sentence, we have $\vdash_\ea \exists x\s(h(x)=i)\to\Box(\exists x\s(h(x)=i))$. Now we see:
\begin{align*}
\vdash_\ea\ &\exists x (h(x)=i)\wedge\Box\left([B]\vee \bigvee_{j\in U_b}\exists y\s (h(y)=j)\right)\\
&\to \Box\left(\exists x\s (h(x)=i)\wedge\left([B]\vee\bigvee_{j\in U_b}(\exists y\s (h(y)=j))\right)\right)\\
&\to \Box\triangle\left([B]\vee \bigvee_{j\in U_{b-1}} \exists y\s(h(y)=j)\right)\\
&\to \Box\left([B]\vee \bigvee_{j\in U_{b-1}} \exists y\s(h(y)=j)\right),
\end{align*}
where the final step holds since $[B]\vee \bigvee_{j\in U_{b-1}} \exists y\s(h(y)=j)$ is equivalent, over $\ea$, to a $\Sigma_1$-sentence.
Now we can apply this repeatedly to \cref{leftovers} in order to obtain
\begin{align*}
\vdash_T \exists x\s (h(x)=i) &\to \Box\left([B]\vee \bigvee_{j\in U_0}\exists y\s (h(y)=j)\right)\\
&\leftrightarrow \Box([B]\ \vee\perp)\\
&\leftrightarrow \Box[B],
\end{align*}
where we used that $U_0=\emptyset$.

Since this holds for all $i\in\lval\Box B\rval$, we can conclude that $\vdash_T [\Box B]\to\Box[B]$, as desired.
\end{prf}

\subsection{The Completeness Theorem}\label{sec_completeness_theorem}

In this section, we formulate and prove our completeness theorem in its abstract form. First, we define provability logics.
\begin{defi}
Let $T$ be a theory and let $P(x)$ be a $\Sigma_1$-formula in one free variable. If $A$ is an $\mathcal{L}$-sentence, we write $\Box A$ for $P(\lg A\rg)$.
\begin{itemize}\index[term]{Realization@($\Sigma_1$-)realization}\index[symb]{asigma@$\sigma$, $\sigma_P$, $\sigma_\Box$, $\sigma_T$, $\sigma^f_T$}\index[term]{Logic@($\Sigma_1$-)logic}\index[term]{Provability logic@(Fast) ($\Sigma_1$-)provability logic}
\item[(i)]	A \emph{realization} is a function $\sigma$ that assigns, to each proposition letter $p$ in $\mathcal{L}_\Box$, an $\mathcal{L}$-sentence $\sigma(p)$. We call $\sigma$ a $\Sigma_1$-realization if $\sigma(p)\in\Sigma_1$ for all $p$.
\item[(ii)]	Given a realization $\sigma$, we define the function $\sigma_P$ from $\mathcal{L}_\Box$ to $\mathcal{L}$-sentences by:
\begin{itemize}
\item[(a)]	$\sigma_P(\bot)$ is $\bot$ and $\sigma_P(p)$ is $\sigma(p)$ for every proposition letter $p$;
\item[(b)]	$\sigma_P(B\circ C)$ is $\sigma_P(B)\circ\sigma_P(C)$ for all $B,C\in\mathcal{L}_\Box$ and $\circ\in\{\wedge,\vee,\to\}$;
\item[(c)]	$\sigma_P(\Box B)$ is $\Box(\sigma_P(B))$ for all $B\in\mathcal{L}_\Box$.
\end{itemize}
\item[(iii)]	The \emph{logic for $P$} w.r.t.\@ $T$ is defined as the set of all $A\in\mathcal{L}_\Box$ such that $\vdash_T \sigma_P(A)$ for every realization $\sigma$. The $\Sigma_1$-logic for $P$ w.r.t.\@ $T$ is the set of all $A\in\mathcal{L}_\Box$ such that $\vdash_T \sigma_P(A)$ for every $\Sigma_1$-realization $\sigma$.
\end{itemize}
By abuse of notation, we also write $\sigma_\Box$ instead of $\sigma_P$, and we say `logic for $\Box$' instead of `logic for $P$'. We write $\sigma_T$ for $\sigma_{\Box_T} = \sigma_{\bew_T}$ and $\sigma^f_T$ for $\sigma_{\Box^f_T} = \sigma_{\bew^f_T}$. The ($\Sigma_1$-)logic for $\Box_T$ w.r.t.\@ $T$ is called the \emph{\textup{(}$\Sigma_1$-\textup{)}provability logic of $T$}, and the ($\Sigma_1$-)logic for $\Box^f_T$ w.r.t.\@ $T$ is called the \emph{fast \textup{(}$\Sigma_1$-\textup{)}provability logic of $T$}. \ruitje
\end{defi}

Now, we again adopt the conventions and notation from \cref{sec_Solovay_function}. All the work from \cref{sec_preservation} now leads to the following result.
\begin{thm}\label{embedding}
Define the $\Sigma_1$-realization $\sigma$ by $\sigma(p)=[p]$ for every proposition letter $p$. Then $\vdash_T \sigma_\Box(A)\leftrightarrow[A]$ for all $A\in\mathcal{L}_\Box$.
\end{thm}
\begin{prf}
This follows by induction on the complexity of $A$ using \cref{vee_pres}, \cref{wedge_pres}, \cref{to_pres} and \cref{box_pres}.
\end{prf}
The following result tells us what the `real' behaviour of the Solovay function $h$ is, in the case that $T$ is $\Sigma_1$-sound.
\begin{prop}\label{jammer_harry}
Suppose that $T$ is $\Sigma_1$-sound. Then $\mathbb{N}\models h(x)=0$.
\end{prop}
\begin{prf}
Since $\prec$ is conversely well-founded, we know that $h$ must have a certain limit $i\in\mathbb{N}$. Then $\exists x\s(h(x)=i)$ is a true $\Sigma_1$-sentence, which means that $\vdash_T\exists x\s(h(x)=i)$. Now suppose that $i>0$. Then $\vdash_T \exists x\s (h(x)=i)\to \Box(\exists y\s (i\int h(y)))$ by \cref{danger_of_leaving_home}, so we must have $\vdash_T\Box(\exists y\s (i\int h(y)))$. Since $\Box(\exists y\s (i\int h(y)))$ is a $\Sigma_1$-sentence and $T$ is $\Sigma_1$-sound, we see that $\mathbb{N}\models\Box(\exists y\s (i\int h(y)))$. By requirement (ii) for a good pair, we get $\vdash_T \exists y\s (i\prec h(y))$. Since $\exists y\s (i\prec h(y))$ is a $\Sigma_1$-sentence and $T$ is $\Sigma_1$-sound, we get $\mathbb{N}\models \exists y\s (i\prec h(y))$. However, this is impossible as $i$ is supposed to be the limit of $h$. So $i=0$, and the result follows.
\end{prf}
Now we can finally formulate and prove our main result.
\begin{thm}\label{Completeness_Theorem}
Let $T$ be a $\Sigma_1$-sound theory. Suppose we have a good pair $(P,Q)$ for $T$ such that $\vdash_T\textup{CP}_P$. Then the \textup{(}$\Sigma_1$-\textup{)}logic for $P$ is equal to the set of theorems of $\iGLC$.
\end{thm}
\begin{prf}
As before, let us abbreviate $P(\lg A\rg)$ as $\Box A$, for $\mathcal{L}$-sentences $A$. Since $P(x)$ is a provability predicate for $T$, we see that the ($\Sigma_1$-)logic for $\Box$ contains the axioms of $\mathsf{iGL}$ and is closed under $\to$E and Nec. Since $\vdash_T\text{CP}_\Box$, we see that the ($\Sigma_1$-)logic for $\Box$ also contains all sentences of the form $A\to\Box A$, where $A\in\mathcal{L}_\Box$. So the ($\Sigma_1$-)logic for $\Box$ contains all theorems of $\iGLC$.

Now suppose that we have $A\in\mathcal{L}_\Box$ such that $\mathsf{iGLC}\nvdash A$. Then by \cref{completeness_iGLC}, there exists a finite, irreflexive, realistic model $M_0 = \langle W_0,\int_0,\r_0,V_0\rangle$ in which $A$ is not valid. We label the nodes of $M_0$ as $W_0=\{1, \ldots, r\}$ in such a way that $M_0,r\nVdash A$. By shrinking $W_0$ to $\{i\in W_0\mid r\int_0 i\}$ if necessary, we may assume without loss of generality that $r$ is the $\int_0$-least element of $W_0$.

Now define the model $M$, the Solovay function $h$, and the $\Sigma_1$-sentences $[B]$ for $B\in\mathcal{L}_\Box$ as above. It is easy to show that $M_0,i\Vdash B$ iff $M,i\Vdash B$ for all $B\in\mathcal{L}_\Box$ and all $i$ with $1\leq i\leq r$. So we have $M,r\nVdash A$, that is, $r\not\in \lval A\rval$. Now define the $\Sigma_1$-realization $\sigma$ by $\sigma(p)=[p]$ for every proposition letter $p$. By \cref{embedding}, we have $\vdash_T \sigma_\Box(B)\leftrightarrow[B]$ for all $B\in\mathcal{L}_\Box$.

Now suppose for the sake of contradiction that $\vdash_T\sigma_\Box(A)$. Then we also get $\vdash_T[A]$. Since $[A]$ is (equivalent to) a $\Sigma_1$-sentence and $T$ is $\Sigma_1$-sound, we see that $\mathbb{N}\models [A]$. By \cref{jammer_harry}, we also know that $\mathbb{N}\models h(x)=0$. This implies that $0\in \lval A\rval$. However, we also have $0\int r$ and $r\not\in\lval A\rval$, which yields a contradiction. We conclude that $A$ is not in the ($\Sigma_1$-)logic for $\Box$, as desired.
\end{prf}

\section{Applications of the Completeness Theorem}\label{chap_applications}

In the previous section, we proved a completeness theorem in a very abstract form. In this section, we provide several applications of this theorem. In particular, we will determine the fast provability logics of the theories $U^\ast$, for $\Sigma_1$-sound theories $U$, and we will determine the fast and ordinary $\Sigma_1$-provability logics of $\ha$. First of all, we lay some further groundwork in \cref{sec_NNIL}. Then, in \cref{sec_fast_provlog}, we determine the fast provability logics mentioned above. Finally, in \cref{sec_main_result}, we determine the $\Sigma_1$-provability logic of $\ha$.

\subsection{The Sets $\nnil$ and $\tnnil$}\label{sec_NNIL}

In the sequel, $\mathcal{L}_p$ is the language of propositional logic, and for $A\in\mathcal{L}_p$, we write `$\vdash_\mathsf{IPC} A$' to indicate that $A$ is provable in intuitionistic propositional logic.\index[term]{Language!of intuitionistic propositional logic}\index[symb]{alp@$\mathcal{L}_p$}\index[symb]{avdashipc@$\vdash_\mathsf{IPC}$} We notice that, if $\sigma$ is a substitution, $A\in\mathcal{L}_p$ and $P(x)$ is $\Sigma_1$-formula, then $\sigma_P(A)$ does not actually depend on $P$. So we will just write $\sigma(A)$ instead of $\sigma_P(A)$. We will also drop the brackets in expressions of the form $\sigma(A)$ and $\sigma_T(A)$.

Like the authors of \cite{ArdeshirMojtahedi}, we introduce the set of $\nnil$-sentences. 
\begin{defi}\index[symb]{annil@$\nnil$}
The set $\nnil\subseteq\mathcal{L}_p$ (`no nested implications on the left') is defined recursively, as follows:
\begin{itemize}
\item[(i)]	all proposition letters are in $\nnil$, as is $\bot$;
\item[(ii)]	if $A,B\in\nnil$, then $A\wedge B, A\vee B \in\nnil$;
\item[(iii)]	if $A\in\mathcal{L}_p$ contains no implications and $B\in\nnil$, then $A\to B\in\nnil$. \ruitje
\end{itemize}
\end{defi}
That is, a $\nnil$-sentence is a propositional sentence in which no implication occurs in the antcedent of another implication. In the paper \cite{NNIL}, we find the following result, that we will not prove here.
\begin{thm}\label{NNIL_algorithm}\index[symb]{aat3@$A^\ast$}\index[term]{nnilalgorithm@$\nnil$-algorithm}
There exists a computable function $(\cdot)^\ast\colon \mathcal{L}_p\to\nnil$, called the \emph{$\nnil$-algorithm}, such that for every $A\in\mathcal{L}_p$, the following hold:
\begin{itemize}
\item[\textup{(}i\textup{)}]	$\vdash_\mathsf{IPC} A^\ast\to A$;
\item[\textup{(}ii\textup{)}]	if $B\in\nnil$ and $\vdash_\mathsf{IPC} B\to A$, then $\vdash_\mathsf{IPC} B\to A^\ast$;
\item[\textup{(}iii\textup{)}]	if $\sigma$ is a $\Sigma_1$-realization, then $\vdash_\ha \Box_\ha(\sigma A)\leftrightarrow\Box_\ha(\sigma A^\ast)$.
\end{itemize}
\end{thm}
\begin{rem}
Consider the preorder $(\mathcal{L}_p,\leq)$, where $\leq$ is defined by: $A\leq B$ if and only if $\vdash_\mathsf{IPC} A\to B$, for $A,B\in\mathcal{L}_p$. Consider also the subpreorder $(\nnil,\leq)$. Then items (i) and (ii) above say that the $\nnil$-algorithm is left adjoint to the inclusion $\nnil\to\mathcal{L}_p$. \ruitje
\end{rem}
We can get an analogue of (iii) for fast provability.
\begin{cor}\label{NNIL_algorithm_fast}
Let $A\in\mathcal{L}_p$ and let $\sigma$ be a $\Sigma_1$-realization. Then 
\[\vdash_\ha \Box_\ha^f(\sigma A)\leftrightarrow \Box_\ha^f(\sigma A^\ast).\]
\end{cor}
\begin{prf}
Since $\bew^f_\ha$ is a provability predicate for $\ha$, we can derive from \cref{NNIL_algorithm}(iii) that
\[\vdash_\ha \Box^f_\ha(\sigma A) \leftrightarrow \Box^f_\ha\Box_\ha(\sigma A) \leftrightarrow \Box^f_\ha\Box_\ha(\sigma A^\ast) \leftrightarrow \Box^f_\ha(\sigma A^\ast),\]
where we also used \cref{facts_iterated_box}(vii).
\end{prf}

The $\nnil$-algorithm behaves nicely with respect to the theories $U^T$ and $\Sigma_1$-realizations.
\begin{prop}\label{nice}
Suppose that $U$ and $T$ are theories such that $\vdash_\ha \bew_U(x)\to\Box_U\bew_T(x)$. Then for all $\Sigma_1$-realizations $\sigma$ and $C\in\nnil$, we have 
\begin{align*}
\vdash_\ha\Box_{U^T}(\sigma C)\to\Box_U(\sigma C).
\end{align*}
\end{prop}
\begin{prf}
We notice that, since $\sigma$ is a $\Sigma_1$-realization and $C\in\nnil$, we have that $\sigma C$ is equivalent, over $\ea$, to a sentence in $\mathcal{A}$. Using \cref{proofs_under_T} and \cref{preservation}, we see that $\vdash_\ea (\sigma C)^T\to\sigma C$. We get $\vdash_\ea \Box_U(\sigma C)^T\to \Box_U(\sigma C)$. Finally, we notice that the conditions of \cref{formalized_prov_in_U_and_U^T} hold for $V\equiv\ha$, so we get $\vdash_\ha\Box_{U^T} (\sigma C)\to \Box_U(\sigma C)^T \to \Box_U(\sigma C)$, as desired.
\end{prf}

Following \cite{ArdeshirMojtahedi}, we now extend the notion of `no nested implication on the left' to modal sentences.
\begin{defi}\index[symb]{atnnil@$\tnnil$}
The set $\mathsf{TNNIL}\subseteq\mathcal{L}_\Box$ (`thoroughly no nested implications on the left') is defined by recursion, as follows:
\begin{itemize}
\item[(i)]	all proposition letters are in $\mathsf{TNNIL}$, as is $\bot$;
\item[(ii)]	if $A,B\in\mathsf{TNNIL}$, then $A\wedge B, A\vee B, \Box A\in\mathsf{TNNIL}$;
\item[(iii)] if $A,B\in\mathsf{TNNIL}$ and $A$ contains no implications outside the scope of a $\Box$, then also $A\to B\in\mathsf{TNNIL}$. \ruitje
\end{itemize}
\end{defi}
We notice that every $A\in\mathcal{L}_\Box$ can, in a unique way, be written as $C(\vec{p},\Box B_1, \ldots, \Box B_k)$, for certain $C(\vec{p},q_1, \ldots, q_k)\in\mathcal{L}_p$ and \emph{distinct} $B_1, \ldots, B_k\in\mathcal{L}_\Box$. It is easy to show that, with this notation, we have $A\in\mathsf{TNNIL}$ if and only if $C\in\mathsf{NNIL}$ and $B_i\in\mathsf{TNNIL}$ for $1\leq i\leq k$. Now we define an operation on modal formulae as in \cite{ArdeshirMojtahedi}.
\begin{defi}\index[term]{tnnilalgorithm@$\tnnil$-algorithm}\index[symb]{aat4@$A^+$}
The \emph{$\tnnil$-algorithm} $(\cdot)^+\colon\mathcal{L}_\Box\to\mathsf{TNNIL}$ is defined by recursion, as follows. For $A\in\mathcal{L}_\Box$, write $A=C(\vec{p},\Box B_1, \ldots, \Box B_k)$, where $C(\vec{p},q_1, \ldots, q_k)\in\mathcal{L}_p$ and $B_1, \ldots, B_k\in\mathcal{L}_\Box$ are distinct. Then
\[A^+ := C^\ast(\vec{p},\Box B_1^+, \ldots, \Box B_k^+).\wruitje\]
\end{defi}
Notice that, since all the $B_i$ have lower complexity than $A$, the operation $(\cdot)^+$ is well-defined. The following lemmata show how certain results about $\nnil$ and $(\cdot)^\ast$ can be transferred to $\tnnil$ and $(\cdot)^+$. We notice that \cref{handling_plus}(i) also occurs in \cite{ArdeshirMojtahedi} as Corollary 4.7.1.
\begin{lem}\label{handling_plus}
Let $A\in\mathcal{L}_p$ and let $\sigma$ be a $\Sigma_1$-realization. Then the following hold:
\begin{itemize}
\item[\textup{(}i\textup{)}]	$\vdash_\ha \Box_\ha(\sigma_\ha A)\leftrightarrow\Box_\ha(\sigma_\ha A^+)$;
\item[\textup{(}ii\textup{)}]	$\vdash_\ha \Box^f_\ha(\sigma^f_\ha A)\leftrightarrow\Box^f_\ha(\sigma^f_\ha A^+)$.
\end{itemize}
\end{lem}
\begin{prf}
(i)	We proceed by strong induction on the boxdepth of $A$. As above, we write $A$ as $C(\vec{p},\Box B_1, \ldots, \Box B_k)$, where $C(\vec{p},q_1, \ldots, q_k)\in\mathcal{L}_p$ and $B_1, \ldots, B_k\in\mathcal{L}_\Box$ are distinct. Then all the $B_i$ have smaller boxdepth than $A$, so we assume by induction hypothesis that
\begin{align}\label{IH}
\vdash_\ha  \Box_\ha  (\sigma_\ha  B_i) \leftrightarrow \Box_\ha  (\sigma_\ha  B_i^+)\quad\mbox{for}\quad 1\leq i\leq k.
\end{align}
If $\vec{p}=p_1,\ldots,p_l$, then we write $\sigma\vec{p}$ as a shorthand for $\sigma(p_1),\ldots,\sigma(p_l)$. Now we take a $\Sigma_1$-realization $\tau$ such that $\tau\vec{p} = \sigma\vec{p}$ and $\tau(q_i) = \Box_\ha(\sigma_\ha B_i)$ for $1\leq i\leq k$. Now we observe that
\begin{align*}
\sigma_\ha  A &= C(\sigma\vec{p}, \Box_\ha  (\sigma_\ha B_1), \ldots, \Box_\ha  (\sigma_\ha B_k)) = \tau C,\\
\sigma_\ha  A^+ &= C^\ast(\sigma\vec{p}, \Box_\ha  (\sigma_\ha B_1^+), \ldots, \Box_\ha  (\sigma_\ha B_k^+))\quad\mbox{and}\\
\tau C^\ast &= C^\ast(\sigma\vec{p}, \Box_\ha  (\sigma_\ha B_1), \ldots, \Box_\ha  (\sigma_\ha B_k)).
\end{align*}
So \cref{IH} gives $\vdash_\ha  \sigma_\ha  A^+\leftrightarrow \tau C^\ast$. Since $\bew_\ha$ is a provability predicate for $\ha$, we conclude that
\[\vdash_\ha  \Box_\ha  (\sigma_\ha  A)\leftrightarrow\Box_\ha  (\tau C)\leftrightarrow\Box_\ha  (\tau C^\ast)\leftrightarrow\Box_\ha (\sigma_\ha A^+),\]
where we used \cref{NNIL_algorithm}(iii). This completes the induction.

(ii) The proof is completely analogous, but with an appeal to \cref{NNIL_algorithm_fast} instead of \cref{NNIL_algorithm}(iii).
\end{prf}

\begin{lem}\label{handling_TNNIL}
Let $U$ and $T$ be theories such that $\ha\subseteq U\subseteq T$.
\begin{itemize}
\item[\textup{(}i\textup{)}]	Suppose that we have
\begin{align}\label{conservation_realizations_of_NNIL}
\vdash_\mathsf{HA} \Box_U (\sigma C)\leftrightarrow \Box_{U^T} (\sigma C)
\end{align}
for all $C\in\nnil$ and $\Sigma_1$-realizations $\sigma$. Then
\[\vdash_\mathsf{HA} \Box_U (\sigma_U A)\leftrightarrow \Box_{U^T} (\sigma_{U^T} A)\]
for all $A\in\tnnil$ and $\Sigma_1$-realizations $\sigma$.
\item[\textup{(}ii\textup{)}]	Suppose that we have
\begin{align}\label{conservation_realizations_of_NNIL_fast}
\vdash_\mathsf{HA} \Box^f_U (\sigma C)\leftrightarrow \Box^f_{U^T} (\sigma C)
\end{align}
for all $C\in\nnil$ and $\Sigma_1$-realizations $\sigma$. Then
\[\vdash_\mathsf{HA} \Box^f_U (\sigma^f_U A)\leftrightarrow \Box^f_{U^T} (\sigma^f_{U^T} A)\]
for all $A\in\tnnil$ and $\Sigma_1$-realizations $\sigma$.
\end{itemize}
\end{lem}
\begin{prf}
First of all, we observe that the conditions of \cref{HA_and_Ttranslation}(i) are satisfied, so we have $\ha\subseteq U^T$ as well.

(i) We proceed by strong induction on the boxdepth of $A$. Write $A=C(\vec{p},\Box B_1, \ldots, \Box B_k)$, where $C(\vec{p},q_1, \ldots, q_k)\in\mathsf{NNIL}$ and $B_1, \ldots, B_k\in\mathsf{TNNIL}$ are distinct. Then all the $B_i$ have smaller boxdepth than $A$, so we assume by induction hypothesis that
\begin{align}\label{IH_2}
\vdash_\mathsf{HA}\Box_U (\sigma_U B_i)\leftrightarrow \Box_{U^T} (\sigma_{U^T} B_i)\quad\mbox{for}\quad 1\leq i\leq k.
\end{align}
Now we take a $\Sigma_1$-realization $\tau$ such that $\tau\vec{p} = \sigma\vec{p}$ and $\tau(q_i) = \Box_U (\sigma_UB_i)$ for $1\leq i\leq k$. Now we observe that
\begin{align*}
\sigma_U A &= C(\sigma\vec{p}, \Box_U (\sigma_UB_1), \ldots, \Box_U (\sigma_UB_k)) = \tau C\quad\mbox{and}\\
\sigma_{U^T} A &= C(\sigma\vec{p}, \Box_{U^T} (\sigma_{U^T}B_1), \ldots, \Box_{U^T} (\sigma_{U^T}B_k)).
\end{align*}
So \cref{IH_2} gives $\vdash_\mathsf{HA} \sigma_{U^T} A\leftrightarrow \tau C$. Since $\ha\subseteq U^T$, we also get $\vdash_{U^T} \sigma_{U^T} A\leftrightarrow \tau C$. We also know that $\bew_{U^T}$ is a provability predicate for $U^T$, so we also find $\vdash_\ha \Box_{U^T} (\sigma_{U^T}A)\leftrightarrow\Box_{U^T}(\tau C)$. Using \cref{conservation_realizations_of_NNIL}, we get
\[\vdash_\mathsf{HA} \Box_U (\sigma_U A)\leftrightarrow\Box_U (\tau C)\leftrightarrow\Box_{U^T} (\tau C)\leftrightarrow\Box_{U^T} (\sigma_{U^T} A),\]
which completes the induction.

(ii) The proof is again completely analogous, but with an appeal to \cref{conservation_realizations_of_NNIL_fast} instead of \cref{conservation_realizations_of_NNIL}.
\end{prf}

\subsection{Some Fast ($\Sigma_1$-)Provability Logics}\label{sec_fast_provlog}

Let $U$ be a $\Sigma_1$-sound theory. By \cref{Sigma_and_U^T}, the theory $U^\ast$ is also $\Sigma_1$-sound. In order to apply the completeness theorem from \cref{sec_completeness_theorem}, we need to prove the following result.
\begin{lem}\label{goodpair_fast}
The pair $\left(\bew^f_{U^\ast}(x), \bew_{U^\ast}(x)\right)$ is good for $U^\ast$.
\end{lem}
\begin{prf}
By \cref{facts_iterated_box}(v), we know that $\bew^f_{U^\ast}$ is a provability predicate for $U^\ast$, and we also know that $\bew_{U^\ast}$ is a provability predicate for $U^\ast$.

Since $U^\ast$ is $\Sigma_1$-sound, we see by \cref{facts_iterated_box}(vi) that $\mathbb{N}\models\Box^f_{U^\ast} A$ implies that $\mathbb{N}\models \Box_{U^\ast}A$, which implies $\vdash_{U^\ast} A$, for all $\mathcal{L}$-sentences $A$.

By \cref{U^ast_proves_CP_U^ast}, we have $\vdash_{U^\ast} \text{CP}_{U^\ast}$.

The final requirement for a good pair follows from \cref{facts_iterated_box}(vii).
\end{prf}
\begin{thm}\label{CT_HA^ast_fast}
Let $U$ be a $\Sigma_1$-sound theory. Then the fast \textup{(}$\Sigma_1$-\textup{)}provability logic of $U^\ast$ is equal to the set of theorems of $\iGLC$.
\end{thm}
\begin{prf}
By \cref{U^ast_proves_CP_U^ast}, we have $\vdash_{U^\ast} A\to\Box_{U^\ast} A\to\Box^f_{U^\ast} A$ for all $\mathcal{L}$-sentences $A$. Now the statement follows from  \cref{Completeness_Theorem} and \cref{goodpair_fast}.
\end{prf}
\begin{rem}\index[symb]{ahastar@$\ha^\ast$}\index[symb]{apastar@$\pa^\ast$}\label{provlog_pa^ast}
Since $\pa$ is a classical theory, we have $\vdash_\pa B\vee(B\to A)$ for all $\mathcal{L}$-formulae $A$ and $B$. This means that we also have
\[\vdash_\pa \Box_\pa A^\pa \to (B^\pa \vee ((B^\pa\to A^\pa)\wedge\Box_\pa(B^\pa\to A^\pa)))\]
for all $\mathcal{L}$-formulae $A$ and $B$. This, in turn, implies that
\[\vdash_{\pa^\ast} \Box_{\pa^\ast} A\to (B\vee (B\to A)),\]
for all $\mathcal{L}$-formulae $A$ and $B$. So the ($\Sigma_1$-)provability logic of $\pa^\ast$ contains at least the theorems of $\iGLC$ extended with the axiom scheme $\Box A\to (B\vee(B\to A))$. This scheme is called the \emph{propositional trace principle}, or PTP for short. The theory $\iGLC+\text{PTP}$ for IML is sound and complete with respect to finite frames $\langle W,\int,\r\rangle$, such that $w\r v$ iff $w\prec v$ for all $w,v\in W$. The first author showed in \cite{completeness} that the ($\Sigma_1$-)provability logic of $\pa^\ast$ contains \emph{exactly} the theorems of $\iGLC+\text{PTP}$. Since $\iGLC+\text{PTP}$ is a proper extension of $\iGLC$, we have an example of a theory for which the fast and ordinary provability logics do not coincide.

Presently, the provability logic for ordinary provability of $\ha^\ast$ is unknown.
As pointed out by Mojtaba Mojtahedi to us in correspondence, it strictly extends {\sf iGLC}.
 A simple example is the principle: 
 \[\Box(\Box \bot \to (\neg\, A \to (B\vee C))) \to \Box (\Box\bot \to ((\neg\, A \to B) \vee (\neg\, A \to C))).\]
Mohammad Ardeshir and Mojtaba Mojtahedi have a manuscript, soon to be published, that gives a characterization of the $\Sigma_1$-provability logic of $\ha^\ast$. 
 \ruitje
\end{rem}

We now turn our attention to determining the fast $\Sigma_1$-provability logic of $\ha$.
\begin{thm}\label{CT_fast_HA}
Let $A\in\mathcal{L}_\Box$. Then $A$ is in the fast $\Sigma_1$-provability logic of $\ha$ if and only if $\iGLC\vdash A^+$.
\end{thm}
\begin{rem}
This result gives an `indirect' characterization of the fast $\Sigma_1$-provability logic of $\ha$, since we first have to apply the $\tnnil$-algorithm, and then see whether the result is provable in $\iGLC$. But we can already see that the fast $\Sigma_1$-provability logic of $\ha$ is decidable, since $\iGLC$ is decidable (this follows from the proof of \cref{completeness_iGLC}). In the paper \cite{ArdeshirMojtahedi}, the authors give a direct characterization of the set $\{A\in\mathcal{L}_\Box\mid\ \vdash_\iGLC A^+\}$, by providing an axiomatization for it. \ruitje
\end{rem}
\begin{prf}[Proof of \cref{CT_fast_HA}]
First of all, we show that the conditions of \cref{handling_TNNIL}(ii) are satisfied with $U\equiv T\equiv\ha$. Let $C\in\nnil$ and let $\sigma$ be a $\Sigma_1$-realization. We observe that the conditions of \cref{HA_and_Ttranslation}(ii) are satisfied with $U\equiv V\equiv T\equiv \ha$, so we have that $\vdash_\ha \Box_\ha A\to \Box_{\ha^\ast} A$ for all $\mathcal{L}$-formulae $A$. We also see that the conditions of \cref{nice} are satisfied with $U\equiv T\equiv \ha$. Combining these two, we see that
\[\vdash_\ha \Box_\ha(\sigma C)\leftrightarrow\Box_{\ha^\ast}(\sigma C).\]
Since $\Box_{\ha^\ast}(\sigma C)$ is a $\Sigma_1$-sentence, we also have $\vdash_\ha \Box_{\ha^\ast}(\sigma C)\leftrightarrow (\Box_{\ha^\ast}(\sigma C))^\ha$. Since $\bew^f_\ha$ is a provability predicate for $\ha$, we get $\vdash_\ha \Box^f_\ha\Box_\ha(\sigma C)\leftrightarrow\Box^f_\ha(\Box_{\ha^\ast}(\sigma C))^\ha$. Using \cref{facts_iterated_box}(vii) and \ref{formalized_prov_in_U^T_and_U_fast}, we get
\begin{align*}
\vdash_\ha \Box^f_\ha (\sigma C) &\leftrightarrow \Box^f_\ha\Box_\ha (\sigma C)\\
&\leftrightarrow \Box^f_\ha(\Box_{\ha^\ast} (\sigma C))^\ha\\
&\leftrightarrow \Box^f_{\ha^\ast}\Box_{\ha^\ast} (\sigma C)\\
&\leftrightarrow \Box^f_{\ha^\ast}(\sigma C),
\end{align*}
as desired.

Using \cref{handling_plus}(ii) and \cref{handling_TNNIL}(ii), we now see that
\[\vdash_\ha \Box^f_\ha(\sigma^f_\ha A) \leftrightarrow \Box^f_\ha(\sigma^f_\ha A^+) \leftrightarrow \Box^f_{\ha^\ast} (\sigma^f_{\ha^\ast} A^+)\]
for all $A\in\mathcal{L}_\Box$ and $\Sigma_1$-realizations $\sigma$. Since $\ha$ is sound, we see that $\mathbb{N}\models \Box^f_\ha(\sigma^f_\ha A)$ if and only if $\mathbb{N}\models \Box^f_{\ha^\ast} (\sigma^f_{\ha^\ast} A^+)$. We also know that $\ha^\ast$ is $\Sigma_1$-sound, so using \cref{facts_iterated_box}(vi), we can now see that
\[
\vdash_\ha \sigma^f_\ha A\quad\mbox{iff}\quad\mathbb{N}\models\Box^f_\ha (\sigma^f_\ha A) \quad\mbox{iff}\quad\mathbb{N}\models \Box^f_{\ha^\ast} (\sigma^f_{\ha^\ast} A^+)\quad\mbox{iff}\quad\vdash_{\ha^\ast} \sigma^f_{\ha^\ast} A^+.\]
This means that $A$ is in the fast $\Sigma_1$-provability logic of $\ha$ if and only if $A^+$ is in the fast $\Sigma_1$-provability logic of $\ha^\ast$. By \cref{CT_HA^ast_fast}, the latter holds if and only if $\vdash_\iGLC A^+$.
\end{prf}

\subsection{A Theory with {\sf iGLC} as Provability Logic}
In this  section, we present an arithmetical theory that has {\sf iGLC} as its provability logic for ordinary provability. 

Recall the theory slow Heyting Arithmetic $\sha$, that satisfies $\sha=\ha$ and $\sha\leq\ha$, but not $\ha\leq\sha$. We consider the theory $\hahat := \ha^\sha$. By \cref{HA_and_Ttranslation}(i), we have $\ha\subseteq\hahat$ and by \cref{Sigma_and_U^T}, the theory $\hahat$ is $\Sigma_1$-sound. Moreover, by \cref{prov_in_hahat}, we know that $\vdash_\ha \Box_\hahat A\leftrightarrow \Box_\ha A^\sha$ for all $\mathcal{L}$-formulae $A$.\index[symb]{ahahat@$\hahat$}

We show that the ($\Sigma_1$-)provability logic of this theory is equal to the set of theorems of $\iGLC$. In order to do this, we need to find a good pair for $\hahat$. In the previous section, the role of $P(x)$ was fulfilled by \emph{fast} provability. In this section, we put ordinary provability for $\hahat$ here. For $Q(x)$, we take $\bew_{\sha^\ast}$. We know from \cref{formalized_prov_in_U_and_U^T} with $U\equiv T\equiv \sha$ and $V\equiv \ea$ that $\vdash_\ea \Box_{\sha^\ast} A\leftrightarrow \Box_\sha A^\sha$ for all $\mathcal{L}$-formulae $A$.
\begin{lem}\label{goodpair}
The pair $\left(\bew_\hahat(x), \bew_{\sha^\ast}(x)\right)$ is good for $\hahat$.
\end{lem}
\begin{prf}
We already know that $\bew_\hahat$ is a provability predicate for $\hahat$. Moreover, since $\sha = \ha$, we also have $\sha^\sha = \ha^\sha$, that is, $\sha^\ast = \hahat$. Since $\bew_{\sha^\ast}(x)$ is a provability predicate for $\sha^\ast$, it must also be a provability predicate for $\hahat$.

Next, let $A$ be an $\mathcal{L}$-sentence. We know from \cref{basic_facts_BoxT}(i) that $\mathbb{N}\models \Box_\hahat A$ implies $\vdash_\hahat A$.

Moreover, by \ref{U^ast_proves_CP_U^ast} with $U\equiv\ha$ and $T\equiv\sha$, we have $\vdash_{\hahat} \text{CP}_{\sha^\ast}$.

Finally, let $S$ be a $\Sigma_1$-sentence. By \cref{Sigma_and_U^T}, we have $\vdash_\ea S\leftrightarrow S^\sha$. We also have that $\vdash_\ea \Box_{\sha^\ast} S \leftrightarrow (\Box_{\sha^\ast} S)^\sha$. Now we use \cref{stealing_from_PA}(i) to find that:
\begin{align*}
\vdash_\ha \Box_\hahat\Box_{\sha^\ast} S &\leftrightarrow \Box_\ha(\Box_{\sha^\ast} S)^\sha\\
&\leftrightarrow \Box_\ha\Box_{\sha^\ast} S\\
&\leftrightarrow \Box_\ha\Box_\sha S^\sha\\
&\leftrightarrow \Box_\ha\Box_\sha S\\
&\to \Box_\ha S\\
&\leftrightarrow \Box_\ha S^\sha\\
&\leftrightarrow \Box_\hahat S.
\end{align*}
Since $\ha\subseteq\hahat$, the final requirement for a good pair follows.
\end{prf}
Now that we have our good pair, we can prove the following.
\begin{thm}\label{CT_hahat}
The \textup{(}$\Sigma_1$-\textup{)}provability logic of $\hahat$ is exactly the set of theorems of $\iGLC$.
\end{thm}
\begin{prf}
Since $\sha\leq \ha$, and $\vdash_\hahat \text{CP}_{\sha^\ast}$, we see that 
\[\vdash_\hahat A\to \Box_{\sha^\ast} A \to \Box_\sha A^\sha \to \Box_\ha A^\sha\to \Box_\hahat A\] for every $\mathcal{L}$-sentence $A$. This means that $\vdash_\hahat\text{CP}_\hahat$, so both statements follow from \cref{Completeness_Theorem} and \cref{goodpair}.
\end{prf}

\subsection{The $\Sigma_1$-Provability Logic of $\ha$}\label{sec_main_result}
We use Theorem~\ref{CT_hahat} to
 determine the (ordinary) $\Sigma_1$-provability logic of $\ha$. This  is 
  the main result of the paper \cite{ArdeshirMojtahedi}, but the authors arrive at it using different methods.

\begin{thm}\label{CT_HA}
Let $A\in\mathcal{L}_\Box$. Then $A$ is in the $\Sigma_1$-provability logic of $\ha$ if and only if $\iGLC\vdash A^+$.
\end{thm}
\begin{prf}
First, we show that the conditions of \ref{handling_TNNIL}(i) are satisfied with $U\equiv \ha$ and $T\equiv\sha$. Let $C\in\nnil$ and let $\sigma$ be a $\Sigma_1$-realization. By \cref{HA_and_Ttranslation}(ii) with $U\equiv V\equiv \ha$ and $T\equiv\sha$, we have that $\vdash_\ha \Box_\ha A\to\Box_\hahat A$ for all $\mathcal{L}$-formulae $A$. We also see that the conditions of \ref{nice} are satisfied with $U\equiv\ha$ and $T\equiv \sha$, so we see that
\[\vdash_\ha \Box_\ha(\sigma C)\leftrightarrow\Box_\hahat(\sigma C),\]
as desired.

Now we can use \cref{handling_plus}(i) and \cref{handling_TNNIL}(i) to see that
\[\vdash_\ha \Box_\ha(\sigma_\ha A)\leftrightarrow \Box_\ha(\sigma_\ha A^+) \leftrightarrow \Box_{\hahat}(\sigma_\hahat A^+)\]
for all $A\in\mathcal{L}_\Box$ and $\Sigma_1$-realizations $\sigma$. Since $\ha$ is sound, we get
\[\vdash_\ha\sigma_\ha A\quad\mbox{iff}\quad \mathbb{N}\models\Box_\ha\sigma_\ha A \quad\mbox{iff}\quad\mathbb{N}\models \Box_\hahat(\sigma_\hahat A^+)\quad\mbox{iff}\quad \vdash_\hahat\sigma_\hahat A^+.\]
This means that $A$ is in the $\Sigma_1$-provability logic of $\ha$ if and only if $A^+$ is in the $\Sigma_1$-provability logic of $\hahat$. By \cref{CT_hahat}, the latter holds if and only if $\vdash_\iGLC A^+$.
\end{prf}

\section{Conclusion}

In this paper, our goal was to give a Solovay-style embedding of frames equipped with both an intuitionistic relation $\int $ and a modal relation $\r$. In order to approach this task, we considered theories that prove their own completeness principle. This project has led to the following results and insights.

\begin{itemize}
\item[(i)]	We were able to give a Solovay-style embedding of finite, irreflexive, realistic frames for IML, in the presence of the completeness principle and the principle $\Box\triangle S\to\Box S$ for $S\in \Sigma_1$.
\item[(ii)]	We reproved the result from \cite{ArdeshirMojtahedi} that the $\Sigma_1$-provability logic of Heyting Arithmetic is equal to the set $\{A\in\mathcal{L}_\Box\mid\ \vdash_\iGLC A^+\}$.
\item[(iii)]	We showed that the fast $\Sigma_1$-provability logic of $\ha$ is also equal to this set.
\item[(iv)]	We showed that for any $\Sigma_1$-sound theory $U$, the fast ($\Sigma_1$-)provability logic of $U^\ast$ is equal to the set of theorems of $\iGLC$.
\item[(v)]	We found an intuitionistic theory of arithmetic other than $\pa^\ast$, namely the theory $\hahat$, for which we were able to determine the provability logic, to wit {\sf iGLC}.
\item[(vi)]	We discovered that for the theory $\pa^\ast$, the fast provability logic and the ordinary provability logic do not coincide.
\end{itemize}

\end{document}